\documentclass[12pt,reqno]{amsart}

\usepackage{latexsym}
\usepackage{graphicx}

\newcommand{\tspace}{\mbox{\rule[-7pt]{0pt}{20pt}}}
\newcommand{\T}{{\mathbf T}^m}

\newcommand{\sm}{\setminus}
\newcommand{\szego}{Szeg\"o }

\newcommand{\Si}{\Sigma}
\newcommand{\inv}{^{-1}}
\newcommand{\kahler}{K\"ahler }

\newcommand{\wt}{\widetilde}
\newcommand{\wh}{\widehat}
\newcommand{\PP}{{\mathbb P}}
\newcommand{\N}{{\mathbb N}}
\newcommand{\R}{{\mathbb R}}
\newcommand{\C}{{\mathbb C}}

\newcommand{\Z}{{\mathbb Z}}

\newcommand{\CP}{\C\PP}
\renewcommand{\d}{\partial}
\newcommand{\dbar}{\bar\partial}
\newcommand{\ddbar}{\partial\dbar}

\newcommand{\half}{{\frac{1}{2}}}
\newcommand{\vol}{{\operatorname{Vol}}}

\newcommand{\codim}{{\operatorname{codim\,}}}

\newcommand{\FS}{{{\operatorname{FS}}}}

\renewcommand{\phi}{\varphi}

\newcommand{\ispa}[1]{\langle \,#1 \,\rangle } %%%%%%%%%%%%% Tatsuya added
\newcommand{\tf}[2]{{\textstyle \frac{#1}{#2}}} %%%%%%%%%%%%% Tatsuya added
\newcommand{\comp}{{\scriptstyle \circ}} %%%%%%%%%%%%%%%%%%% Tatsuya added

\newcommand{\ccal}{\mathcal{C}}

\newcommand{\hcal}{\mathcal{H}}
\newcommand{\ical}{\mathcal{I}}

\newcommand{\lcal}{\mathcal{L}}
\newcommand{\mcal}{\mathcal{M}}

\newcommand{\ocal}{\mathcal{O}}

\newcommand{\al}{\alpha}
\newcommand{\be}{\beta}
\newcommand{\ga}{\gamma}
\newcommand{\Ga}{\Gamma}

\newcommand{\la}{\lambda}

\newcommand{\om}{\omega}

\newtheorem{theo}{{\sc Theorem}}[section]

\newtheorem{cor}[theo]{{\sc Corollary}}

\newtheorem{lem}[theo]{{\sc Lemma}}
\newtheorem{prop}[theo]{{\sc Proposition}}

\newenvironment{rem}{\medskip\noindent{\it Remark:\/} }{\medskip}
\newtheorem{defin}[theo]{{\sc Definition}}

\title[Distribution laws for  integrable eigenfunctions]
{Distribution laws for integrable eigenfunctions }

\author{Bernard Shiffman}
\address{Department of Mathematics, Johns Hopkins University, Baltimore, MD
21218, USA} \email{shiffman@math.jhu.edu}

\author{Tatsuya Tate}
\address{Department of Mathematics, Keio University, Keio University
3-14-1 Hiyoshi Kohoku-ku, Yokohama, 223--8522 Japan}
\email{tate@math.keio.ac.jp}

\author{Steve Zelditch}
\address{Department of Mathematics, Johns Hopkins University, Baltimore, MD
21218, USA}
\email{zelditch@math.jhu.edu}

\thanks{Research partially supported by NSF grants DMS-0100474 (first author)
and  DMS-0071358 (third author) and by JSPS (second author).}

\date{June 11, 2003}

\begin{document}

\begin{abstract} We
determine  the asymptotics  of the joint eigenfunctions of the
torus action on a toric K\"ahler variety. Such varieties are
models of completely integrable systems in complex geometry.  We
first determine the pointwise asymptotics of the eigenfunctions,
which show that they behave like  Gaussians centered at the
corresponding classical torus. We then
 show that there is a  universal Gaussian scaling limit of
the distribution function near its center. We also determine the
limit distribution for the tails of the eigenfunctions  on large
length scales. These are not universal but depend on the global
geometry of the toric variety and in particular on the details of
the exponential decay of the eigenfunctions away from the
classically allowed set.

\end{abstract}

\maketitle

\tableofcontents

\section{Introduction} A problem of considerable interest in both mathematics and physics is
to determine  the asymptotics of the   distribution functions
$$D_j (t): = \vol\{z: |\phi_{j} (z)|^2 > t\}$$
of an orthonormal basis $\{\phi_j\}$ of eigenfunctions of a
Laplacian (or similar Hamiltonian) on a compact manifold (see
\cite{Y}). In general, it is hopelessly difficult to obtain more
than crude bounds on such distribution functions, which of course
control the $\lcal^p$-norms of the eigenfunctions.   Numerical
analyses and
 heuristics from quantum chaos and disordered systems
  suggest however a rich picture in which
the asymptotics of $D_j(t)$ is related to the classical dynamics
underlying the eigenvalue problem. (Some references will be
discussed at the end of the introduction.) The purpose of this
paper is to give a rather complete analysis of the limit
distribution of eigenfunctions in one of the few settings where
such a detailed analysis is possible,   namely where the phase
space is a    toric \kahler variety $(M, \omega)$.

Let us recall the definitions (see \S 1 for details). Toric
varieties are complex manifolds on which the complex torus
$(\C^*)^m$ acts with an open dense orbit.  By a toric \kahler
variety we mean a toric variety equipped with a  \kahler form $(M,
\omega)$ that is invariant under the underlying real torus ${\bf
T}^m$.  The action of ${\bf T}^m$ is Hamiltonian with respect to
$\omega$, and thus  toric varieties are models of completely
integrable systems. They are of a very special type because
integrable systems usually generate an $\R^m$ action rather than a
${\bf T}^m$ action. Although there are rigidity theorems limiting
the class of such examples in the world of real Riemannian
manifolds \cite{LS}, toric varieties provide a plentiful
collection in the world of complex manifolds.

The torus action can be `quantized' or linearized on the Hilbert
space completion of the coordinate ring
\begin{equation} {\mathcal H} := \bigoplus_{N=0}^\infty H^0(M, L^N),
\end{equation}
where $L \to M$ is a holomorphic line bundle with $c_1(L) =
\frac{1}{2\pi} \omega$ and where $H^0(M, L^N)$ denotes the space
of holomorphic sections of its $N$-th tensor power. This
quantization is generally known as the holomorphic (Bargmann-Fock)
representation in the physics literature.  The space ${\mathcal H}$  is spanned by joint
  eigenfunctions of the linearized $(\C^*)^m$ action, which we refer to
  as `monomials.' In the fundamental
case of $M = \CP^m$, the joint eigenfunctions are  the monomials
given in an affine chart by \begin{equation} \chi_{\alpha} : \C^m
\to \C,\;\;\; \chi_{\alpha}(z) = z^{\alpha}. \end{equation} The
monomials lift (by homogenization) to homogeneous monomials on
$\C^{m + 1}$.

We consider the case where $M$ is a smooth projective toric
variety; i.e., $M=M_P$, $L=L_P$, where $P$ is an integral Delzant
polytope (see \S\ref{Background}). Then  the linearized ${\bf
T}^m$ action is generated by  $m$ commuting operators $\hat{I}_j$,
$j = 1, \dots, m$ on $M_P$ which preserve $H^0(M_P, L_P^N)$, and
the joint spectrum of the eigenvalue problem
\begin{equation} \label{QCI} \hat{I}_j \phi^P_{\alpha} = \alpha_j \phi^P_{\alpha},\;\;
\alpha=(\al_1,\dots\al_m) \in \R^m,\;\;\; \phi^P_{\alpha}\in H^0(M_P,
L_P^N)
\end{equation}  consists of lattice points $\alpha \in N P \cap \Z^m$.

Our main results concern the asymptotics of their distribution functions
\begin{equation}
\label{DISTRIBUTION}
D_{\gamma}(t):=  \vol\{z\in M_P : |\phi^P_{\gamma}(z)|^2 > t\}
\end{equation}
with $\gamma \in NP \cap \Z^{m}$ as $N \to \infty$. The function $
|\phi^P_{\gamma}(z)|^2$ is often called the `Husimi distribution'
in the physics literature, and thus our results determine its
distribution law. The norm $|\phi^P_{\gamma}(z)|$ of
$\phi^P_{\gamma}(z)\in L_P^N$ is  the pull-back of the
Fubini-Study norm under  a monomial embedding of the form
\begin{equation} \label{PHIP} \Phi_{P}^{c}=[c_{\alpha(1)}\chi_{\alpha(1)},\ldots,c_{\alpha(d+1)}\chi_{\alpha(d+1)}]:
(\C^*)^m \to \CP^{d}\;,\qquad
P\cap\Z^m=\{\al(1),\dots,\al(d+1)\}\;, \end{equation} for a choice
of constants $c_{\al(j)}\in\C^*$.
 The volume in $M_P$ and the Hermitian norm $h_P^c$  on $L_P$ are by definition the  pull-backs
   of the Fubini-Study metric and
form under  this monomial embedding, and the $\lcal^2$ norm on the
space $H^0(M,L^N_P)$ is in turn induced from the volume form and
the Hermitian pointwise norm of $h_P^c$.   (See \S
\ref{Background} for details.)

\begin{figure}[htb]\label{monfig}
\begin{center}
\includegraphics*{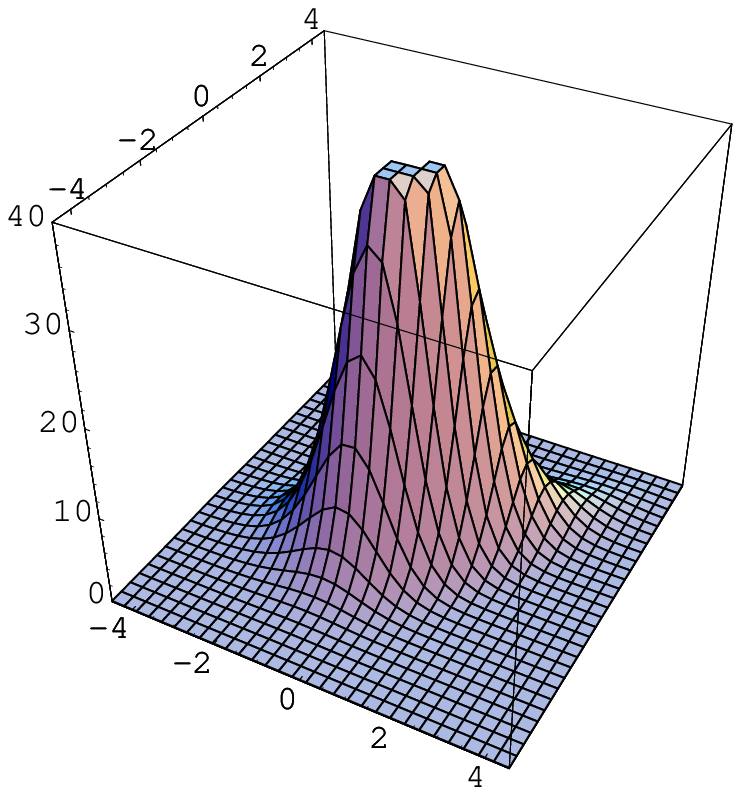}
\hspace{10pt}
\includegraphics*{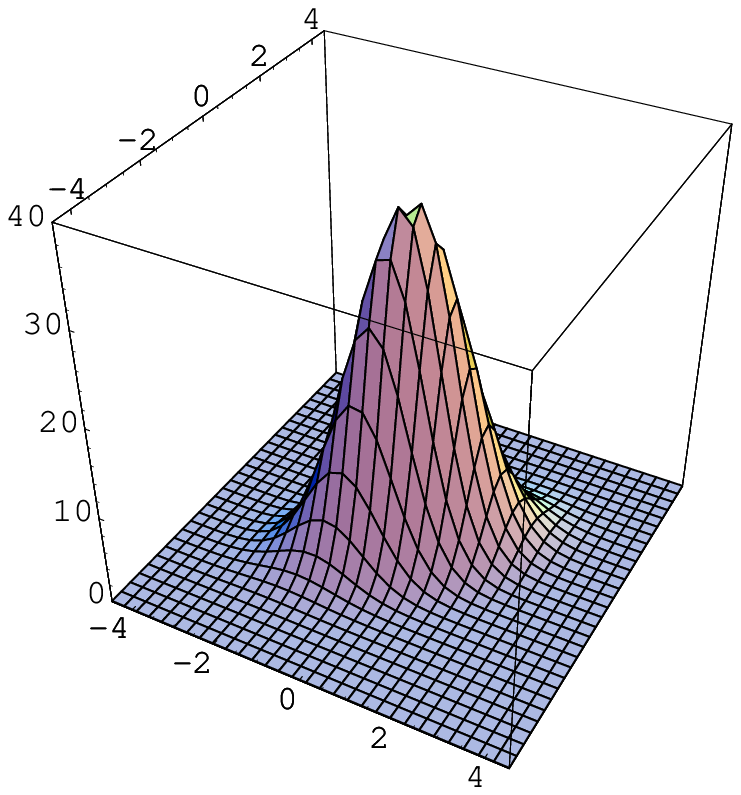}
 \caption{ $2\pi /N$ times the
monomial $|\phi_{N\alpha}^{7\Sigma}(z)|^{2}$ for $N=1$ (left) and $N =
\infty$ (right) for  $m=2$ and $\alpha=(2,3)$,
where $\Sigma$ is the standard simplex. The variable $z$ is chosen
as $u \mapsto z=e^{(\rho_{\alpha}+u/\sqrt{N})/2}$. See Proposition \ref{LOCALP}
in Section \ref{massasym}.
}
\end{center}
\end{figure}

As Figure \ref{monfig} illustrates, the monomial $\phi^P_{\alpha}$ is
something like a   Gaussian bump centered on the real torus
$\mu_P\inv(\al)$, where $\mu_P:M_P\to P$ is the moment map for the
classical Hamiltonian ${\bf T}^m$-action on $M_P$ (see
\S 1).

Since we would like to determine the properties of
eigenfunctions $\phi_\al\in NP$ when $\al$ is large, but not
necessarily a multiple of a lattice point in $P$, we  shall
consider a sequence of approximate multiples, as in  the following
definition:

\begin{defin} Let $\alpha_{N} \in NP
\cap \Z^{m}$ be a sequence of lattice points, and let $x\in P$. We
say that $\{\alpha_{N}\}$ is a   sequence of approximate multiples
of $x$ if
\begin{equation}
\label{genL} \alpha_{N}=Nx +O(1).
\end{equation} \end{defin}

Our first result gives the  pointwise behavior of the
eigenfunctions:

\begin{theo}
Let $x$ be a point in the interior $P^{o}$ of the polytope $P$,
which is not necessarily a lattice point. Then there exists a
non-negative function $b_{x}^{P}
\in\ccal^\infty((\mathbb{C}^{*})^{m})$ such that $b_{x}^P(z)=0$ if
and only if $z \in \mu_{P}^{-1}(x)$, and for every sequence
$\alpha_{N} \in NP \cap \Z^{m}$  of approximate multiples of $x$,
we have
\[
|\phi^{P}_{\alpha_{N}}(z)|^{2}=c(P,x) \left( \frac{N}{2\pi}
\right)^{m/2}{e^{-Nb_{x}^{P}(z)}}\left[1 +O(N^{-1})\right]
\]
uniformly on  $(\mathbb{C}^{*})^{m}$, where $c(P,x)\in\R^+$.
\label{pointww}
\end{theo}

\bigskip

The constant $c(P,x)$ is defined in \eqref{matrixa} and
\eqref{CONSTX} (and also appears in Theorems \ref{unrescal} and
\ref{mldist}).

When we consider a sequence of the lattice points of the form $\alpha_{N}=N\alpha$
with a lattice point $\alpha$,
the assumption that $\alpha \in P^{o}$ is not necessary.
In fact, we will give the pointwise asymptotics for $\alpha_{N}=N\alpha$ with
any lattice point $\alpha \in P \cap \Z^{m}$ in Section \ref{massasym}
(see Propositions \ref{LOCALSS} and \ref{LOCALV}).

In the above Theorem, the function $b_{x}^{P}(z)$ ($x \in P^{o}$)
on $(\mathbb{C}^{*})^{m}$ is defined by
\begin{equation}\label{b}
b_{x}^{P}(z):= \log \left(\frac{ \sum_{\beta \in
P}|c_\be z^\be|^{2}}{ \sum_{\beta \in
P}e^{-\ispa{\tau_{x}^P(z),\beta}}|c_\be z^\be|^{2}}
\right) -\ispa{\tau_{x}^P(z),x},
\end{equation}
where $\tau_{x}^P(z) \in \mathbb{R}^{m}$ is the vector
given by the equation
\begin{equation}\label{tau}
\mu_{P}(e^{-\tau_{x}^P(z)/2}z)=x,\quad z \in
(\mathbb{C}^{*})^{m}\;.
\end{equation}
Here $\mu_{P}$ is the moment map for the $\T$ action on $M_{P}$
(see \eqref{muP} for the definition), and we write
\begin{equation} e^rz=(e^{r_1}z_1,\dots,e^{r_m}z_m) \qquad \mbox{for }\ r\in\R^m,\ z\in
(\C^*)^m\;.
\end{equation}
We can express (as in \cite[(17)--(18)]{SZ2}) the function
$b_{x}^{P}$ in the more intuitive form as follows:  we introduce
the real power `monomials'
$$\quad |\chi_x(z)|:= |z|^x = |z_1|^{x_1} \cdots |z_m|^{x_m} $$ and define
\begin{equation}\label{mcal} \quad \mcal_x^P (z):=\frac
{|\chi_x(z)|_P} {\sup |\chi_x(z)|_P}\;,
\end{equation}
where (cf.\ \eqref{mp})
$$
|\chi_x(z)|_P := \frac{|\chi_x(z)|}{\sqrt{\sum_{\beta \in P}|c_\be
z^\be|^{2}}}\; \;\;\;({z\in (\C^*)^m}).$$ (The normalized monomial
$\mcal^P_x$ has sup-norm $1$, attained on the torus
$\mu_P\inv(x)$.) Then \eqref{b} is equivalent to:
\begin{equation} \label{interpretb} b_x^P(z)=-2\log
{\mcal^P_{q(z)}(z)}\;,\end{equation}

Since the sequence of
monomials flattens out exponentially quickly away from the
peak set $\mu_P\inv (x)$, the distribution function is clearly
tending to zero.  The rate of decay of the distribution function is given by the
following result:

\begin{theo}
\begin{itemize}

\item[{\rm (i)}] Let   $\alpha_{N} \in NP
\cap \Z^{m}$ be a sequence of lattice points which are approximate
multiples of $x \in P^{o}$ (see   \eqref{genL}). Then, for $t>0$,
we have
\[
D_{\alpha_{N}}(t) \sim
\frac{(\pi m)^{m/2}}{c(P,x)\Gamma(m/2 +1)}
\left(
\frac{\log N}{N}
\right)^{m/2}.
\]

\item[{\rm (ii)}] Let $\alpha \in \partial P \cap \Z^{m}$. Then, for $t>0$, we have
\[
D_{N\alpha}(t) \sim \frac {(\pi
d(\alpha))^{d(\alpha)/2}}{{c(P,\alpha)\Gamma(d(\alpha)/2 +1)}}\,
\left(\frac {\log N}{N}\right)^{d(\alpha)/2},
\]
where  $d(\alpha):= m + \codim F_\al$,  $F_\al$ being the face of
$P$ containing $\al$.

\end{itemize}

\label{unrescal}
\end{theo}

\noindent Here $\sim$ means the ratio of the left and right hand
sides tends to $1$, and the constants $c(P,x),\ c(P,\alpha)$ are
given by  \eqref{matrixa} and \eqref{CONSTC}--\eqref{CONSTX}.
We recall that if $x$ is a point in a face $F$ of codimension $r$,
then $\mu_{P}^{-1}(x) \cong {\bf T}^{m-r}$. Hence $d(\alpha)=2m
-\dim \mu_{P}^{-1}(\alpha)$.

 The exponentially
localized behavior of the monomials suggests studying the
distribution function on various length scales. First, we   show
that the $D_{\alpha_{N}}$ have a universal scaling limit on a
small length scale:

\begin{theo}\label{mldist}

\begin{itemize}

\item[{\rm (i)}] Let $\alpha_{N} \in NP \cap \Z^{m}$ be a sequence of lattice points satisfying
the condition \eqref{genL} with some point $x \in P^{o}$.
Then, for $0< t \leq c(P,x)$, we have
\[
\lim_{N \to \infty} (N/2\pi)^{m/2}D_{\alpha_{N}}
\left(
(N/2\pi)^{m/2}t
\right)=
\frac{1}{c(P,x)\Gamma(m/2 +1)}
(\log (c(P,x)/t))^{m/2}.
\]

\item[{\rm (ii)}]Let $\al \in P \cap \Z^{m}$.
Then
\[
\lim_{N \to \infty}\left(
N/{2\pi}\right)^{d(\alpha)/2}D_{N\alpha} \left( \left(
N/{2\pi}\right)^{d(\alpha)/2}t \right)=
\frac{1}{c(P,\alpha)\Gamma(d(\alpha)/2 +1)}
\left(
\log(c(P,\alpha)/t)
\right)^{d(\alpha)/2},
\] for $0<t\le c(P,\al)$, where $d(\alpha)=\codim \mu_{P}^{-1}(\alpha)$.

\end{itemize}

\end{theo}

\begin{figure}[htb]
\centerline{\includegraphics*{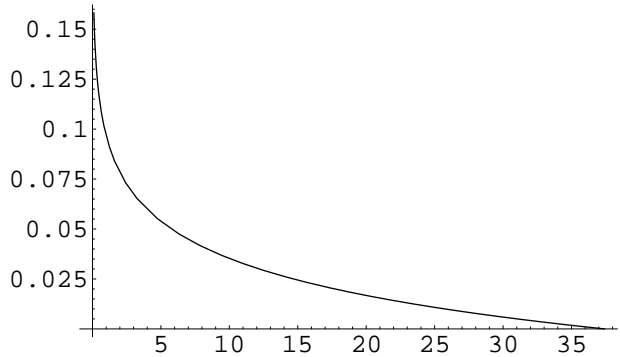}} %scale=0.7
\caption{ Scaling limit distribution for $|\varphi_{N\alpha}^{P}|^{2}$
with $m=2$, $\alpha =(2,3)$, $P=7\Sigma$} \label{limitD}
\end{figure}

A sample graph of the scaling limit distribution function for $P=7\Si$  is given in
Figure
\ref{limitD}.

As in Theorem \ref{mldist}, the limit of the rescaled
distributions has a universal form, {\it i.e.} it does not depend
on  the geometry of the manifold $M_{P}$, and is given by   a
logarithmic power of the form $(\log c/t)^{d/2}$ with some
constant $d$. The  logarithmic power appears because the
$\lcal^{2}$-normalized monomials are close to  Gaussian around the
peak set $\mu_{P}^{-1}(x)$ on a vector space of dimension $m/2$
(or dimension $d(\alpha)$ for the boundary lattice cases). More
precisely,  Theorem \ref{pointww} (and Propositions \ref{LOCALSS}
and \ref{LOCALV}) shows that the  function $b_{x}^{P}$ has a
positive definite Hessian at the peak point. We then observe that
the distribution  of a Gaussian function on $\R^{d}$,
\[
g(u):=\frac{e^{-\ispa{Au,u}/2}}{\sqrt{\det A}},\quad u \in \R^{d},
\]
(where $A$ is a real positive $(d \times d)$-matrix) is given by
the logarithmic power law
\begin{gather*}
\nu_{A}(u \in \R^{d}\,;\,g(u)>t)=
\frac{1}{c\Gamma(d/2 +1)}
\left(
\log \frac{c}{t}
\right)^{d/2}\\[5pt]
c=\frac{1}{\sqrt{\det A}},\quad 0 < t \leq c,
\end{gather*}
relative to  the normalized Lebesgue measure
\[
\nu_{A}:=\frac{\det A}{(2\pi)^{d/2}}\,du.
\]
Thus, the  rescaled distribution of an $\lcal^{2}$-normalized
monomial  at the `center' of its localized  bump has a
 universal Gaussian form.

To analyze the `tails' of the
eigenfunctions,  we next use an exponential rescaling of  the distribution function so
that the global distribution law has a non-zero limit as $N \to \infty$.  As may be
expected, it is no longer universal but depends on the geometry of $(M_{P},
\omega_{P})$.

\begin{theo} Let $\alpha_{N} \in NP \cap \Z^{m}$
satisfy the condition \eqref{genL} with a point $x \in P^{o}$.
Then
\begin{equation*}
\lim_{N \to \infty} D_{\alpha_{N}}(e^{-Nt})= \int_{\{\rho \in
\mathbb{R}^{m}\,;\,b_{x}^{P}(\rho) < t \}} \det
A(\rho)\,d\rho\;, \end{equation*}
where $A(\rho)$ is the Hessian matrix of $\log \sum_{\be\in P}
|c_\be|^2e^{\ispa{\be,\rho}}$.
\label{tnonrescaledintro}
\end{theo}

\noindent A more general scaling limit law is given in Theorem
\ref{tnonrescaled}. In the above theorem, the assumption that $x
\in P^{o}$ is not necessary. In fact, we will give similar result
for $\alpha_{N}=N\alpha$ with any lattice point $\alpha \in P$
(Theorem \ref{thnonrescaled}) in Section \ref{invmprob}.

Our strategy for proving  Theorems \ref{unrescal} and
\ref{tnonrescaledintro} on the  distribution functions of
monomials is based on their  pointwise asymptotics (Theorem
\ref{pointww}, and also Propositions \ref{LOCALSS} and
\ref{LOCALV} in Section \ref{massasym}). Pointwise asymptotics of
monomials are more or less equivalent to asymptotics of their
$L^{2k}$ norms. Since the latter are of independent interest, we
state these asymptotics explicitly in:

\begin{theo}\label{ml2knorm}

\begin{itemize}

\item[{\rm (i)}] Let $\alpha_{N} \in NP \cap \Z^{m}$ satisfy the condition \eqref{genL}
for a point $x \in P^{o}$. Let $\|\varphi_{\gamma}^{P}\|_{2k}$ denote the $\lcal^{2k}$-norm
of the $\lcal^{2}$-normalized monomial $\varphi_{\gamma}^{P}$ with the weight $\gamma \in NP$.
Then we have
\[
\|\varphi_{\alpha_{N}}\|_{2k}^{2k}=
\frac{c(P,x)^{k-1}}{k^{m/2}}
\left(
\frac{N}{2\pi}
\right)^{(k-1)m/2}(1+O_{k}(N^{-1})),
\]
where $O_{k}(N^{-1})$ depends on $k$.

\item[{\rm (ii)}] Let $\alpha \in P \cap \Z^{m}$ with $d(\alpha)=\codim \mu_{P}^{-1}(\alpha)$.
Then we have
\[
\|\wh{\varphi}_{N\alpha}^{P}\|_{2k}^{2k}
=\frac{c(P,\alpha)^{k-1}}{k^{d(\alpha)/2}}
\left(
\frac{N}{2\pi}
\right)^{(k-1)d(\alpha)/2}(1+O_{k}(N^{-1})).
\]

\end{itemize}

\end{theo}

We close the introduction with some general remarks and
references.  As mentioned above,  our results pertain to phase
space distribution of eigenfunctions (Husimi distributions) rather
than to their configuration space distribution. To our knowledge,
the only prior  example in dimension $> 1$ for which the limit
distribution of eigenfunctions has been determined is  the case of
certain (so-called) Hecke eigenfunctions of discrete quantum cat
maps, due to Kurlberg-Rudnick \cite{KR}.
They work in a simpler discrete model rather than
the holomorphic model.  The main result of Kurlberg-Rudnick
\cite{KR} is that the distribution functions of (un-scaled)
eigenfunctions tend to the semi-circle law. Their method was to
relate the eigenfunctions to exponential sums studied by Katz
\cite{Ka} and to apply the value distribution of exponential sums.
 Since value
distribution depends on the representation, it is not  clear that
the same semi-circle law would hold for  Hecke eigenfunctions  in
the holomorphic (Bargmann-Fock) representation, and this appears
to be a challenging and interesting problem. Part of the
motivation for this paper was to set a baseline for eigenfunction
distribution problems by studying a class of explicitly solvable
examples.

It would also be interesting to  study the limit distribution  of
real eigenfunctions in the Schr\"odinger representation, i.e. on
the configuration space rather than the phase space. To our
knowledge, the physics results mainly pertain to these
configuration space results.
 The cases most studied and speculated about are those of chaotic or
disordered systems. When eigenfunctions are  delocalized,  their
spatial distribution is conjectured to be Gaussian (see e.g.
\cite{B, FE, H, HR,M,  SS, PA}). In the opposite regime where the
  eigenfunctions of disordered systems are exponentially localized, the
expected distribution in the low amplitude (tail) region is given
by a power of a logarithm \cite{MF}, precisely the one we obtained
in the high amplitude (center) region. The reason is that the
distribution law  is universal and Gaussian in the tail region for
exponentially localized eigenfunctions of disordered systems
\cite{MF}, while it is universal and Gaussian  in the center for
our problem (Theorem \ref{pointww}).

The only studies we have located which are related to distribution
laws of  integrable eigenfunctions are those of one of the authors
with J. A. Toth (cf. \cite{TZ}) and that  of Berry-Hannay-Ozorio
de Almeida \cite{BHO}, which describes the asymptotic expansions
of $\lcal^{2p}$ norms (moment intensities) of real  oscillatory
integrals of several stable types. Such oscillatory integrals
define quasimodes for a quantum integrable system, and in generic
cases one can express eigenfunctions as oscillatory integrals of
various kinds \cite{TZ}. In general, many possible kinds of
oscillatory integrals could arise, including ones associated to
singular Lagrangean tori. Hence, one can only expect complete
results in special cases. To take the simplest example, our
methods could be adapted to find the scaling limit distributions
of squares of the standard spherical harmonics
 $Y^N_m$ on $S^2$ (or $S^m$). These eigenfunctions are also joint
 eigenfunctions of commuting operators which generate a quantum
 torus action. To our knowledge, the distribution laws of even
 such simple eigenfunctions are unknown at this time.

\vspace{5pt}

\noindent{{\bf Acknowledgments:}} This paper was written during a
stay of the second author at Johns Hopkins University on
a JSPS fellowship. He would
like to express his special thanks to  the faculty in the
Department of Mathematics of Johns Hopkins University.

\section{Background on toric varieties and moment polytopes}

\label{Background} We summarize here some basic facts and
terminology on toric varieties from \cite{STZ}. Recall that a {\it
toric variety\/} is a complex algebraic variety $M$ containing the
complex torus $$(\C^{*})^{m}:= (\C\sm\{0\})\times\cdots\times
(\C\sm\{0\})$$ as a Zariski-dense open set such that the group
action of $(\C^*)^m$ on itself extends to $M$.
We consider here smooth projective toric varieties; they can be
given the structure of a symplectic manifold such that the
restriction of the action to the underlying real torus
$$\T=\{(\zeta_1,\dots\,\zeta_m)\in (\C^*)^m: |\zeta_j|=1, 1\le j\le m\}$$ is a
Hamiltonian action (see \S \ref{s-torus}). These toric varieties
can be constructed from Delzant polytopes either by symplectic
reduction (see \cite{Gu}) or by gluing affine toric varieties
described by the normal fan of the polytope (see \cite{Fu}).
However, for our analysis, it is more convenient to define the
toric variety $M_{P}$ associated to the Delzant polytope $P$
through a {\it monomial embedding\/} as follows (see
\cite[Chapter~5]{GKZ}). Suppose that $P$ is a Delzant polytope and
let
$$P\cap\Z^m=\{\al(1),\al(2),\dots\al(d+1)\}\;.$$
We shall write $\sum_{\al\in P}=\sum_{\al\in P\cap\Z^m}$. (Recall
that a {\it Delzant polytope\/} is a convex integral polytope in
$\R^m$ with the property that each vertex is incident to exactly
$m$ edges and the primitive vectors in $\Z^m$ parallel to these
edges generate $\Z^m$.)

To define the monomial embedding,  we fix an arbitrary
$c=(c_{\alpha(1),\ldots,\alpha(d+1)})
\in (\C^{*})^{d+1}$.  Then, we define
\begin{equation}
\Phi_{P}^{c}=[c_{\alpha(1)}\chi_{\alpha(1)},\ldots,c_{\alpha(d+1)}\chi_{\alpha(d+1)}]:
(\C^*)^m \to \CP^{d}\,; \label{PHPC}
\end{equation}
i.e.,
\[
\Phi_{P}^{c}(z)=[c_{\alpha(1)}z^{\alpha(1)},\ldots,c_{\alpha(d+1)}z^{\alpha(d+1)}],
\quad z \in (\C^*)^m.
\]
The toric variety $M_{P}^{c}=M_{P}$ is defined as the
Zariski-closure of the image $\Phi_{P}^{c}((\C^*)^m)$ of the
monomial embedding $\Phi_{P}^{c}$ in the complex projective space
$\CP^{d}$.

Since our polytope $P$ is assumed to be Delzant, $\Phi_{P}^{c}$ is
an embedding and the variety $M_{P}$ is smooth. The symplectic (or
K\"ahler) form   on $M_P$ is given by
\begin{equation} \om_P^c=\om_P:=
\Phi_P^{c*}\om_\FS\;,\label{omega}\end{equation} where $\om_\FS=
\frac{i}{2\pi}\ddbar \log\|\zeta\|^2$ denotes the Fubini-Study
\kahler form on $\CP^d$ with homogeneous coordinates
$(\zeta_0,\dots,\zeta_d)$.  On $(\C^*)^m$, we have
\begin{equation} \om_P^c=\frac{\sqrt{-1}}{2\pi}\partial
\bar{\partial} \log \sum_{\alpha \in P} |c_{\alpha}|^2
|z^{\alpha}|^2. \label{pullfubini}
\end{equation}
The volume form on $M_P$ is given by $d\vol_M=\frac 1
{m!}\om_P^m$.

\subsection{The line bundle $L_P^c$ and associated circle bundle $X_P^c$}

We define the line bundle $L_{P}^{c} \to M_{P}^{c}$ by
$L_{P}=L_{P}^{c}:=\Phi_{P}^{c*}\ocal(1)$, where $\ocal(1)$ denotes
the hyperplane section bundle on $\CP^d$.  Recall that the space
of holomorphic sections $H^0(\CP^{d},\ocal(1))$ consists of the
linear functions $\la:\C^{d+1}\to\C$, and that the {\it
Fubini-study metric\/}  on $\ocal(1)$ is given by
$$|\la|_\FS ([\zeta])=\frac{|\la(\zeta)| }{\|\zeta\|}\qquad
(\zeta\in\C^{d+1})\;,$$ which has curvature form
$\om_\FS=\frac{i}{2\pi}\ddbar \log \|\zeta\|^2$. We  endow $L_P^c$
with the Hermitian metric $h_P^c:= \Phi_P^{c*} h_{\FS}$, which has curvature $\omega_P^c$.

Each monomial $\chi_{\alpha}$  with $\alpha \in P\cap\Z^m$
corresponds to a section of $H^0(M_P^c, L_P^c)$ and vice versa. To
explicitly define this correspondence, we make the
identifications:
\begin{equation}\label{identify}\chi^{P}_{\al(j)}\equiv
c_{\al(j)}\inv\Phi_P^{c*}\zeta_j\in H^0(M_P^c, L_P^c) \cong
\Phi_{P}^{c*}H^{0}(\CP^{d},\ocal(1)) \;, \quad 1\le j \le
d+1\;.\end{equation} More generally,
\begin{equation}\label{identifyN}\Phi_P^{c*}
H^{0}(\CP^{d},\ocal(N)) \cong H^{0}(M_{P},L_{P}^{N})\qquad (N\ge
1),\end{equation}
 and a basis for $H^{0}(M_{P},L_{P}^{N})$ is given by
the sections $\{\chi^P_\ga:\ga\in NP\cap \Z^m\}$ corresponding to
the monomials $\{\chi_\ga\}$.  These sections are given by
\begin{equation*}
\chi^{P}_{\gamma}=\chi^{P}_{\beta_{1}} \otimes \cdots \otimes
\chi^{P}_{\beta_{N}},
\end{equation*}
where $\beta_{1},\ldots ,\beta_{N} \in P\cap\Z^m$ such that
$\gamma =\beta_{1}+ \cdots +\beta_{N}$ (see\cite{Fu, STZ}).

So far, we have not specified the constants $c_\al$.  For studying
our phenomena, the choice of constants defining the toric variety
$M_P$ is not important. However,  when our polytope $P$ is the
full simplex $p\Si$, we shall use the special choice $c_\al =
{p \choose \al}^{1/2}$, where ${p\choose \al}$ is the multinomial
coefficient (see \S \ref{s-proj}).

The associated principal $S^{1}$-bundle $X_{P}^{c}=X_{P}$ of
the line bundle $L_{P} \to M_{P}$ is defined by
\[
X_{P}^{c}:=\{
(z,v) \in L_{P}^{-1}\,;\,|v|_P =1\},
\]
where $|v|_P$ denotes the norm of $v$ with respect to the
Hermitian metric on $L_{P}^{-1}$ induced by $h_{P}^{c}$.

We identify sections $s_N$ of $L^N$ with equivariant functions
$\wh{s}_N$ on $X$ by the rule
\begin{equation} \label{sNhat}\wh{s}_N(\lambda) = \left( \lambda^{\otimes
N}, s_N(z) \right)\,,\quad \la\in X_z\,.\end{equation} Clearly,
$\wh{s}_N(e^{i \theta} \cdot  x) = e^{iN\theta} \wh{s}_N(x)$ if
$s_N \in H^0(M_P^c, L_P^N)$. It should be noted that for each $s
\in H^{0}(M_{P}^c,L_{P}^{N})$, we have
\[
|\wh{s}_{N}(x)|=|s_N(z)|_{P}
\]
where $x \in X_{P}$ is in the fiber over $z \in M_{P}^c$,
$|s_N(z)|_{P}$ denotes the norm with respect to  the Hermitian
metric on $L_P^N$ induced by the metric $h_{P}^{c}$.

In particular, for $\al\in P\cap \Z^m$,  the `monomial'
$\chi_{\al}^{P}\in H^0(M_P^c, L_P^c)$ given by (\ref{identify})
lifts to an equivariant function $\wh\chi_\al^P$ on the circle
bundle $X_P^c\to M_P^c$, and we write
\begin{equation}\label{mhat}\wh m_{\al(j)}^P:= c_{\al(j)} \wh
\chi_{\al(j)}^P= \zeta_j\circ\iota_P\end{equation} where $\iota_P
:X_P^c \to S^{2d +1}$ is the lift of the  embedding  $M_P^c
\hookrightarrow \CP^{d}$ ($d=\#P-1$), that is, $\iota_{P}$ is the
restriction to $X_{P}$ of the natural inclusion $L_{P}^{-1}
\hookrightarrow \ocal(-1)$. (Of course, $\wh m_\al^P$ depends on
$c$, which we omit to simplify notation.) We also consider the
monomials
$$m_{\al}^{P}:=c_\al \chi_{\al}^{P}$$
so that $\wh m_\al^P$ is the equivariant lift of $m_\al$ to
$X_P^c$.  In terms of local coordinates $(z,\theta)$ on
$\pi\inv((\C^*)^m)\subset X_P^c$, we have
\begin{equation}\label{equivariantm}
\wh m_\al^P(z,\theta)= \frac{e^{i\theta}c_\al
z^\al}{\left(\sum_{\be\in P}|c_\be z^\be|^2\right)^{1/2}}\;.
\end{equation}

 Since its absolute value is
independent of $\theta$, we can write
\begin{equation} \label{mp} |m^P_\al(z)|_P = |\wh m^P_\al(z)|= \frac {|c_\al
z^\al|} {\left(\sum_{\be\in P}|c_\be
z^\be|^2\right)^{1/2}}\;.\end{equation}

We give $H^0(M_P,L_P^N)$ the inner product
\begin{equation}\langle s_1,\bar s_2 \rangle = \int_M \left\langle
s_1(z),\overline{s_2
(z)}\right\rangle_{h_N}d\vol_M(z)\;,\quad\quad s_1,s_2\in
H^{0}(M_{P},L_{P}^{N})\label{inner}\;,\end{equation} and the
$\lcal^2$ norm $\|s\|=\langle s,\bar s \rangle$. We
note that the sections
$$\left\{ \chi_{\alpha}^P\in H^0(M_P, L_P^N) :
\alpha=(\al_1,\dots,\al_m)\in NP\right\}$$ are orthogonal but not
normalized. We normalize them to obtain an orthonormal basis for
$H^0(M_P,L_P^N)$ consisting of the sections
\begin{equation}\label{phi}\phi_{\alpha}^{P} :=\frac{
\chi_{\alpha}^{P}}{\|\chi_{\alpha}^{P}\|}\;.\end{equation}
Their equivariant lifts $\wh\phi_{\alpha}^{P}$ form an
orthonormal set of monomials on $X_P^c$.

We include  below a table of notation to help the reader keep track of the various
monomials:

\bigskip\begin{center}
\begin{tabular}{|c|c|c|c|} \hline \tspace& monomials on $\C^m$ & sections of $L_P^N$ &
monomials on $X_P$\\  \hline  \hline \tspace $\al\in NP$&
$\chi_{\alpha} (z) = z^{\alpha}$ &
$\chi^P_\alpha$   & $\wh\chi^P_\alpha$ \\
 \hline \tspace $\al\in P$& &
$m_{\alpha}^{P} = c_{\alpha} \chi_{\alpha}^{P}$ \quad ($N=1$) & $\wh m_{\alpha}^{P}$ \\
\hline \tspace $\al\in NP$&& $\phi_{\alpha}^{P} =
\chi_{\alpha}^{P}/\|\chi_{\alpha}^{P}\|$  & $\wh \phi_{\alpha}^{P} $\\
\hline

\end{tabular}
\end{center}

\bigskip

\subsection{Moment maps and torus actions}\label{s-torus}

The group $(\C^*)^m$ acts on $M_P^c$ and the subgroup $\T $ acts
in a Hamiltonian fashion. Let us recall the formula for its moment
map $\mu_{P}^c: M_P^c \to \R^m$:  on the open orbit $(\C^*)^m$, we
have
\begin{equation}\label{muP} \mu_{P}(z)=\mu_P^c (z) =
\frac{1}{\sum_{\alpha \in P} |c_{\alpha}|^2 |z^{\alpha}|^{2}}
\sum_{\alpha \in P} |c_{\alpha}|^2 |z^{\alpha}|^{2}\alpha  =
\sum_{\alpha \in P} |\wh m_\al^P(z)|^2\alpha\;. \end{equation} For
any $c$, the image of $M^c_P$ under $\mu_P^c$ equals $P$. The
moment map $\mu_{P}$ is invariant under the ${\bf T}^{m}$-action
on $M_{P}$. By the identification $(\C^{*})^{m} \cong \T \times
\mathbb{R}^{m}$, the moment map $\mu_{P}$ defines a map from
$\mathbb{R}^{m}$ to $P$, which  is a diffeomorphism between
$\mathbb{R}^{m}$ and the interior $P^{o}$ of $P$.

The action of the real torus ${\bf T}^m$ lifts from $M_P^c$ to
$X_P^c$ and combines with the $S^1$ action to define a ${\mathbf
T}^{m+1}$ action on $X_P^c$.  Recall that under the monomial
embedding
$$\Phi_P^c:(\C^*)^m\hookrightarrow M_P^c\hookrightarrow \CP^{d},\qquad
z\mapsto \big[c_{\al(1)}z^{\al(1)},\dots
,c_{\al(d+1)}z^{\al(d+1)}\big]\;,$$ the $\T$ action on $M_P^c
\subset \CP^{d}$ is given by
\begin{equation}\label{actionM}e^{i\phi} \cdot
[\zeta_1,\dots,\zeta_{d+1}]=
\big[e^{i\langle\al(1),\phi\rangle}\zeta_1,\dots,
e^{i\langle\al(d+1),\phi\rangle}\zeta_{d+1}\big]\;.\end{equation}
The action (\ref{actionM}) lifts to an action  on $L_P\inv$:

\begin{equation}\label{actionX}e^{i\phi} \cdot \zeta
= \big( e^{i\langle\al(1),\phi\rangle}\zeta_1,\dots,
e^{i\langle\al(d+1),\phi\rangle}\zeta_{d+1}\big) \;.\end{equation}
Since the circle bundle $X_P^c\subset S^{2d+1}$ is invariant under
this action, (\ref{actionX}) also gives a lift of the action
(\ref{actionM}) to $X_P^c$.

We also have the standard circle action on $X_P^c$:
\begin{equation}\label{circle} e^{i\theta}\cdot \zeta =  e^{i\theta}
\zeta\;,\end{equation} which commutes with the $\T$-action
(\ref{actionX}). Combining (\ref{actionX}) and (\ref{circle}), we
then obtain a ${\mathbf T}^{m+1}$-action on $X_P^c$:
\begin{equation}\label{actionX+}(e^{i\theta}, e^{i\phi_1},\dots, e^{i\phi_m})\bullet
\zeta =  e^{i\theta}(e^{i\phi}\cdot \zeta)\;.\end{equation}

\subsection{Fourier analysis} \label{FAnalysis}

In this section, we shall explain an aspect of Fourier analysis on toric varieties
which describe the complete integrability of the system.

The Hardy space
$\hcal^2(X_P^c)$ is the Hilbert space spanned by the equivariant lifts of sections to
functions on $X_P^c$ with the inner
product
$$\langle f, \bar g \rangle = \int_X f\bar g dV,\;\;\; dV = \alpha_P \wedge
(d\alpha_P)^{n-1}, $$ where $\alpha_P$ is a contact 1-form defined by
the Hermitian connection on $L_{P}\inv$ such that $d\alpha_P=\pi^{*}\omega_{P}$.
Under the identification $\hcal\cong \hcal^2(X_P^c)$, the inner
product is the same as the inner product on $\hcal = \bigoplus_N H^0(M_P,L_P^N)$ given
by (\ref{inner}). Alternately, $\hcal^2(X_P^c)$ consists of the functions $F\in
\lcal^2(X_P^c)$ satisfying $\dbar_b F=0$ (see e.g., \cite{SZ1,Z}).

Under the
$S^1$ action, the Hardy space then has the orthogonal
decomposition
\begin{equation} \hcal^2(X_P^c)=
\bigoplus_{N=0}^{\infty}\hcal^2_N(X_P^c)\;, \end{equation}
where $\hcal^2_N(X_P^c)$ consists of  elements $\hat s\in \hcal^2(X_P^c)$  such that
$\hat s(e^{i\theta}\cdot x)= e^{iN\theta} \hat s(x)$. We recall that the (equivariant)
`\szego projectors'
$\Pi_{N}$ are the orthogonal projection onto $H^0(M_P^c, L_P^{c N}) \cong
\hcal_{N}^{2}(X_{P}^{c})$. If $\{S^N_j\}$ denotes an orthonormal basis of
$H^{0}(M_{P}^{c},L_{P}^{cN})$, and $\wh{S}_{j}^{N}$ denote their lifts to $X$,  then
the projector $\Pi_{N}$ is given by the kernel function
\begin{equation}\label{szego}\Pi_{N}(x,y) = \sum_{j = 1}^{k_N}
\wh{S}_{j}^{N}(z) \overline{\wh{S}_{j}^{N}(y)}\;: \lcal^2(X_P^c) \to
\hcal^2_N(X_P^c).\end{equation}

We now describe how one can combine the  eigenvalue problems given by \eqref{QCI} for
varying
$N$ into a homogeneous scalar eigenvalue problem on $X$. To do this,  we define
for each $N \in \mathbb{N}$, the `homogenization'
$\wh{NP} \subset \mathbb{Z}^{m+1}$
of the lattice point in the polytope $NP$ to be the set of all lattice point
$\wh{\alpha}^{N}$ of the form
\[
\wh{\alpha}^{N}=\wh{\alpha}:=(Np -|\alpha|,
\alpha_{1},\ldots,\alpha_{m}),\quad
\alpha=(\alpha_{1},\ldots,\alpha_{m}) \in NP \cap \mathbb{Z}^{m},
\]
where $p \geq \max_{\beta \in P \cap \mathbb{Z}^{m}}|\beta|$.
We also define the cone $\Lambda_{P}=\bigcup_{N=1}^{\infty}\wh{NP}$.
It is well known that rays $\N \wh{\alpha}$ in this cone define
a semiclassical limit.

In this section, we use the more precise notation
$\wh{\phi}_{\hat{\alpha}}^P (x)$ for the $\lcal^2$-normalized monomial
$\wh{\phi}_{{\alpha}}^P (x)$ (since $N$ is not specified in the latter), for
$\wh{\alpha} \in \Lambda_{P}$.

The torus action on $X_P^c$ can be quantized to define an action
of the torus as unitary operators on $\hcal^2(X_P^c).$
Specifically, we let $\hat I_1,\dots\hat I_m$ denote the differential
operators on $X_P^c$ generated by the $\T$ action:
\begin{equation}\label{Xi}(\hat I_j \hat S)
(\zeta)=\frac{1}{i}\frac{\d}{\d\phi_j}  \hat
S(e^{i\phi}\cdot\zeta)|_{\phi=0}\;,\quad \hat S\in
\ccal^\infty(X_P^c)\;.
\end{equation}

We recall the following observation from \cite{STZ}:
\begin{prop} For  $1\le j\le m$,

\begin{itemize}

\item[\rm (i)] $\ \hat I_j:\hcal^2_N(X_P^c) \to \hcal^2_N(X_P^c)$;

\item[\rm (ii)] The lifted monomials $\wh{\phi}_{\hat{\alpha}}^P  \in \hcal^2_N(X_P^c)$
satisfy
$\hat I_j \wh{\phi}_{\hat{\alpha}}^P  = \alpha_j \wh{\phi}_{\hat{\alpha}}^P $
($\wh\al\in \Lambda_P$).

\end{itemize}

\end{prop}

Furthermore, we note that
\begin{equation}\label{dtheta} \frac{\d}{\d\theta} :\hcal^2_N(X_P^c) \to
\hcal^2_N(X_P^c)\;,\qquad \frac{1}{i}\frac{\d}{\d\theta}\hat s_N =
N\hat s_N \quad \mbox{for }\ \hat s_N\in
\hcal^2_N(X_P^c)\;.\end{equation}
Thus, the monomials $\wh{\phi}_{\hat \al}$ are the solutions of the joint eigenvalue
problem
\begin{equation} \label{jointeigen}
\hat{I}_j \wh{\phi}_{\hat \al} = \wh \al_j \wh{\phi}_{\hat \al},\;\;
\wh{\alpha} \in \R^{m + 1} ,\;\;\; \bar{\partial}_b
\wh{\phi}_{\hat \al} = 0, \;\;
j=0,\ldots,m\end{equation}
with the commuting operators:
\begin{equation}
\hat{I}_0=\frac{p}{i}\frac{\partial}{\partial \theta}
-\sum_{j=1}^{m}\hat I_{j},\quad \hat I_1,\ \dots\ ,\ \hat I_m\;. \label{Qtorus}
\end{equation}
The joint eigenvalues are the lattice points $\wh{\alpha} \in \Lambda_{P}$.

\subsection{Monomials on projective space}\label{s-proj}

In the case of $\CP^m$, the $\lcal^q$-norms of the monomials
$\wh{\phi}_{N\alpha}^{p\Si}$ can be evaluated explicitly in an
elementary way. We give the details in this section.

When the polytope $P$ is the unit simplex $\Si$, we have
$M_\Si=\CP^m$; furthermore $L_\Si$ is the hyperplane bundle
$\ocal(1)$. We can identify $L_{\Sigma}^{-1}=\ocal_{\CP^{m}}(-1)$
with $\mathbb{C}^{m+1}$ with the origin blown up, and the circle
bundle $X_\Si \subset L_{\Sigma}^{-1}$ is identified with the unit
sphere $S^{2m +1} \subset \mathbb{C}^{m+1}$. The equivariant lifts
to $X_\Si$ of sections of $\ocal(N)=L_\Si^N$ consist of
homogeneous polynomials
$$F(\zeta_0, \dots,
\zeta_m) = \sum_{|\la| = N} C_\la \zeta^\la \qquad (\zeta^\la =
\zeta_0^{\la_0} \cdots \zeta_m^{\la_m})$$ in $m+1$ variables. The
induced Fubini-Study metric on  $\ocal(N)$ is given by
$$|F(\zeta)|_{\Si}= |F(\zeta)|/\|\zeta\|^N, \quad\mbox{for }\ F\in
H^0(\CP^m,\ocal(N))\;.$$ (Here the constants $c_\al$ are taken to
be $1$.) Identifying $F$ with the polynomial
$f(z)=F(1,z_1,\dots,z_m)$, the  norm can be written
$$|f(z)|_\Si=|f(z)|/(1+\|z\|^2)^{N/2} \qquad (z\in \C^m)\;,$$
and the   inner product is given by:
\begin{equation}\label{IP}\langle f, \bar g \rangle =  \frac{1}{m!}
\int _{\C^m}\frac{\langle
f(z),\overline{g(z)}\rangle}{(1+\|z\|^2)^p} \,\om_\FS^m(z),\quad
f,g\in  H^0(\CP^m, \ocal(p)).\end{equation}

We are interested here in the case $P=p\Si$. Then
$M_{p\Sigma}=\CP^{m}$, and the line bundle $L_{p\Sigma}$ is
identified with the $p$-th tensor power $\ocal(p)$ of the
hyperplane section bundle. The circle bundle $X_{p\Sigma}$ is the
lens space $X_{p\Sigma}=S^{2m+1}/\mathbb{Z}_{p}$. By lifting
equivariant functions from the lens space
$X_{p\Sigma}=S^{2m+1}/\mathbb{Z}_{p}$ to the sphere $S^{2m+1}$, we
see that $\hcal_{N}^{2}(X_{p\Sigma})\cong
\hcal_{Np}^{2}(S^{2m+1})$.  Hence, we shall replace the $N$-th
Hardy space $\hcal_{N}^{2}(X_{p\Sigma})$  by
$\hcal_{Np}^{2}(S^{2m+1})$ below. Then the equivariant lift
$\wh{\chi}_{\alpha}^{p\Si}:S^{2m+1} \to \mathbb{C}$ of
$\chi_{\alpha}^{p\Si} \in H^{0}(\CP^{m},\ocal(p))$ is given by the
homogenization:
\[
\wh{\chi}_{\alpha}^{p\Sigma}(x)=x^{\wh{\alpha}},\quad
\wh{\alpha}=(p-|\alpha|,\alpha_{1},\ldots,\alpha_{m}).
\]

In this case, we shall use the special choice of the coefficients
of the monomial embedding:
$$c^*_\al = {p\choose
\al}^\half\;,\qquad {p\choose\al}:=
\frac{p!}{(p-|\al|)!\alpha_1!\cdots \alpha_m!}\;,$$ so that by
\eqref{pullfubini}, we have
\begin{equation*} \om_{p\Si}^{c^*}=\frac{\sqrt{-1}}{2\pi}\partial
\bar{\partial} \log \left(\textstyle{\sum_{|\alpha| \le p}
{p\choose \al} |z^{\alpha}|^2}\right) =
\frac{\sqrt{-1}}{2\pi}\partial \bar{\partial} \log (1+\|z\|^2)^p =
p \om_\FS.
\end{equation*}
Furthermore,
\begin{equation}\label{special} \mu_{p\Si}(z):=\mu_{p\Si}^{c^*}
(z)=  \frac{1}{\sum_{|\alpha|\le p} {p\choose\al}|z^{\alpha}|^{2
}} \sum_{|\alpha|\le p} {p\choose\al}|z^{\alpha}|^{2 }\alpha
=\frac{p}{1+\sum|z_j|^2}(|z_1|^2,\dots ,|z_m|^2)\,,\end{equation}
where the last equality follows by differentiating the identity
$(1+\sum x_j)^p=\sum_{|\al|\le p} {p\choose\al}x^\al$.  Note that
this choice gives us the scaling formula
$$\mu_{p\Si}=p\mu_\Si.$$

As before, for each $m$-dimensional multi-index $\beta$ with $|\beta| \leq p$,
we define an $(m+1)$-dimensional multi-index $\wh{\beta}$ by
\[
\widehat{\beta}=(\widehat{\beta}_{0},\widehat{\beta}_{1},\ldots,\widehat{\beta}_{m}),\quad
\widehat{\beta}_{0}=p-|\beta|,\quad \widehat{\beta}_{j}=\beta_{j}
\quad (j=1,\ldots,m).
\]
Recall that
$|\wh{\phi}_{\beta}^{p\Si}(x)|=|\phi_{\beta}^{p\Si}(z)|_{p\Sigma}$ with $x \in
S^{2m+1}$, $\pi(x)=z \in \CP^{m}$.

The $\lcal^{q}$-norms of the monomial $\wh{\phi}_{N\alpha}^{p\Si}$
can be evaluated explicitly as follows.

\begin{prop} For $\al\in p\Si\cap\Z^m$, we have the precise formula:
\[
\|\wh{\phi}_{N\hat \al}^{p\Si}\|_{q}^{q}= \left[
\frac{(Np+m)!}{(N\widehat{\alpha})!} \right]^{q/2}
\frac{\prod_{j=0}^m\Gamma(Nq\wh{\al}_j/2+1)}{p^{m(q/2-1)}\Gamma(Npq/2+m+1)}.
\]

\label{exact}
\end{prop}

\begin{proof}

We write $\wh\chi_{N\wh\al}=\wh\chi_{N\hat \al}^{p\Si}$, which we
consider as a function in $\hcal^{2}_{Np}(S^{2m+1})$. For $q\ge
1$, we set
\[
\ical_{q}(N)=\int_{\mathbb{C}^{m+1}}e^{-|x|^{2}}|\widehat{\chi}_{N\hat
\al}(x)|^{q}\,d\ell (x),
\]

where $d\ell(x)$ denotes Lebesgue measure on
$\mathbb{C}^{m+1}$. We shall compute the integral $\ical_{q}(N)$
in two ways. First we use polar coordinates on $\mathbb{C}^{m+1}$.
The measure $d\ell$ is expressed as
\[
d\ell (x) =r^{2m+1}\,dr\,d\sigma, \quad
x=r\sigma,\quad r>0,\quad \sigma \in S^{2m+1}.
\]
Then we have
\begin{eqnarray*}
\ical_{q}(N)&=&\int_{0}^{\infty}e^{-r^{2}} r^{Npq+2m+1}\,dr \
\int_{S^{2m+1}}|\widehat{\chi}_{N\hat \al}(\sigma)|^{q}\,d\sigma\\
&&=\half\Ga(Npq/2+m+1) \int_{S^{2m+1}}|\widehat{\chi}_{N\hat
\al}(\sigma)|^{q}\,d\sigma\,.\end{eqnarray*}

To relate the volume form $d\sigma$ on $S^{2m+1}$ to that on
$X_{p\Si}$, we recall that $\om_{p\Si}=p\om_\FS$ and therefore the
volume form on $X_{p\Si}$ is $p^{m}$ times the Fubini-Study
volume. Recalling our convention that $\vol(X_\Si)=\vol (\CP^m) =
\frac 1{m!}$, we then have $d\vol_{X_{p\Si}}=\frac
{p^m}{2\pi^{m+1}}\,d\sigma$, and hence
$$\ical_{q}(N) =\frac{\pi^{m+1}}{p^m}\Ga(Npq/2+m+1) \|\widehat{\chi}_{N\hat
\al}\|^{q}_q.$$

On the other hand, if we use polar coordinates for each
component $x_{j}$ of $x \in \mathbb{C}^{m+1}$, we have
\[
\begin{split}
\ical_{q}(N)
&=\prod_{j=0}^{m}\int_{\mathbb{C}}e^{-|x|^2}|x|^{N{\alpha}_{j}q}\,d\ell(x) \\
&=(2\pi)^{m+1}\prod_{j=0}^{m}\int_{0}^{\infty}e^{-r^{2}}r^{N{\alpha_{j}}q+1}\,dr
=\pi^{m+1}\prod_{j=0}^m \Gamma(N\al_jq/2+1).
\end{split}
\]
Hence we obtain
\begin{equation}\label{kth}
\|\widehat{\chi}_{N\hat \al}\|_{q}^{q} = p^m\;
\frac{\prod_{j=0}^m\Gamma(Nq\al_j/2+1)}{\Gamma(Npq/2+m+1)}
\end{equation}
and therefore
\begin{equation}\label{l2}
\|\widehat{\chi}_{N\hat \al}\|^{2}
=\frac{p^m(N\widehat{\alpha})!}{(Np+m)!}.
\end{equation}
Since $\wh{\phi}_{N\hat \al}^{p\Si}=\widehat{\chi}_{N\hat
\al}/\|\widehat{\chi}_{N\hat \al}\|$, the identity follows from
\eqref{kth}--\eqref{l2}.
\end{proof}

The following proposition is a direct consequence of Stirling's formula
and Proposition \ref{exact}.

\begin{prop}
\label{Stirling} Let $\alpha \in p\Sigma \cap \Z^{m}$, and
let $J=\{j\in\Z: 0\le j\le m,\ \al_j\neq 0\}$, where
$\al_0=p-|\al|$. Set $d(\alpha)=2m-\#J$. Then we have
\[
\|\wh{\varphi}_{N \alpha }^{p\Sigma}\|_{q}^{2q} \sim \left(
\frac{N}{2\pi} \right)^{(q/2-1)d(\alpha)}
\frac{p^{q/2-1}(2\pi)^{r(q-2)}}
{(q/2)^{d(\alpha)}\left(\prod_{j\in J}\al_j\right)^{q/2-1}};
\]
\[
\|\wh{\varphi}_{N\alpha}^{p\Sigma}\|_{\infty}^{2}= p^{-m}
\frac{(Np+m)!}{(N\wh{\alpha})!} \left[ \frac{\prod_{j\in
J}\al_j^{\al_j}} {p^{p}} \right]^{N}.
\]
In particular, again by Stirling's formula, we have
\[
\|\wh{\varphi}_{N\alpha}^{p\Sigma}\|_{\infty}^{2} \sim (2\pi)^{r}
\left(\frac{p}{\prod_{j\in J}\al_j}\right)^{1/2} \left(
\frac{N}{2\pi} \right)^{d(\alpha)/2}.
\]

\end{prop}

In the next section, we will obtain similar formulas for monomials
on general toric varieties.

\section{Pointwise asymptotics on general toric varieties}

\label{massasym}

We now consider the case of general toric varieties. Our first
purpose is to find the pointwise asymptotics of the monomials and
to prove Theorem \ref{pointww}. We then asymptotically determine
their $\lcal^{2k}$-norms.

\subsection{Pointwise asymptotics: interior points}

First of all, we shall consider a sequence $\alpha_{N}$ of lattice points
in $NP$. We assume that there exists a point $x \in P^{o}$ such that
\begin{equation}
\label{GENL}
\alpha_{N}/N =x +O(N^{-1}).
\end{equation}
Under the assumption \eqref{GENL}, the point $\alpha_{N}/N$ is in the interior $P^{o}$
of the polytope $P$ for sufficiently large $N$.
Thus the analysis is performed on the open orbit $(\C^{*})^{m}$, and hence the coordinate
\[
z=e^{\rho/2 +i\varphi} \in (\C^{*})^{m},\quad
\rho, \varphi \in \R^{m}
\]
will be useful. The moment map $\mu_{P}$ is also invariant under the
Hamiltonian ${\bf T}^{m}$-action, and it is well-known (\cite{Fu}) that
it induces a diffeomorphism:
\begin{equation}
\label{momd}
\bar{\mu}_{P}:\R^{m} =(\C^{*})^{m}/{\bf T}^{m} \to P^{o},\quad
\bar{\mu}_{P}(\rho):=\mu_{P}(e^{\rho/2}).
\end{equation}
In these coordinates, the function $b_{x}^{P}$ defined by
\eqref{b} can be written simply as
\begin{equation}
\label{binv}
b_{x}^{P}(\rho)=f(x,\rho)-f(x,\rho_{x}^{P}),\quad
\rho_{x}^{P}=\bar{\mu}_{P}^{-1}(x),\quad
f(x,\rho):=\log k(\rho)-\ispa{\rho,x},
\end{equation}
where the function $k(\rho)$ is the `polytope character'
\begin{equation}
\label{charP}
k(\rho)=\sum_{\beta \in P \cap \Z^{m}}
|c_{\beta}|^{2}e^{\ispa{\rho,\beta}}.
\end{equation}
The vector $\tau_{x}^{P}(z)$ in \eqref{tau} is given by
\begin{equation}
\label{tau22}
\tau_{x}^{P}(z)=\rho-\rho_{x}^{P},\quad z=e^{\rho/2 +i\theta}.
\end{equation}
We also define the  symmetric real $m \times m$ matrix
\begin{equation}
A(P,x):=\sum_{\beta \in P} \left|\wh m_\be^P(e^{\rho_{x}^P/2})\right|^2\beta
\otimes \beta-x \otimes x.
\label{matrixa}
\end{equation}

\begin{lem}
The  real symmetric matrix
$$
A(\rho):= \sum_{\beta
\in P} |\wh{m}_{\beta}^{P}(e^{\rho/2})|^2 \beta \otimes  \beta -
\mu_{P}(e^{\rho/2}) \otimes \mu_{P}(e^{\rho/2})
$$
is positive definite, for all $z \in (\mathbb{C}^{*})^{m}$.

\label{posmat}
\end{lem}

\begin{proof} We must show that $(\la\otimes\la,A(\rho)) >0$ for $\la\in
(\R^m)'\sm\{0\}$. Consider the vectors $u,v\in\R^{d+1}$  given by $u_\al=|\wh
m_\al^P(e^{\rho/2})|$,
$v_\al=|\wh m_\al^P(e^{\rho/2})|\la(\al)$. Since $\|u\|^2=\sum _{\al\in P} |\wh
m_\al^P(e^{\rho/2})|^2=1$ and $\mu_P(e^{\rho/2})=\sum_{\al\in P}|\wh
m_\al^P(e^{\rho/2})|^2\al$, we have
\begin{eqnarray}\label{S}(\la\otimes\la,A(\rho)) &=& \sum_{\al\in P} |\wh
m_\al^P(e^{\rho/2})|^2 \la(\al)^2 - \left(\sum _{\al\in P} |\wh
m_\al^P(e^{\rho/2})|^2 \la(\al)\right)^2\nonumber \\&=& \|v\|^2-
\langle u, v\rangle^2 \ =\ \|u\|^2\|v\|^2- \langle u, v\rangle^2\
\ge \ 0\;.\end{eqnarray} Since $u_\al=|\wh m_\al^P(e^{\rho/2})|\ne
0$ for all $\al\in P\cap\Z^{m}$ and $v_\al/u_\al=\la(\al)$ is not
constant on $P\cap\Z^{m}$,  the Cauchy-Schwartz inequality in
(\ref{S}) is strict.
\end{proof}

In particular, the matrix
$$A(P,x)=  A(\rho_{x}^{P})$$
is positive definite.

 For a  point $x \in P^{o}$, we
now define the constant
\begin{equation}
\label{CONSTX} c(P,x):=\frac{1}{\sqrt{\det A(P,x)}}.
\end{equation}

\begin{lem}

In the coordinates $z=e^{\rho/2+i\theta}$ on $(\mathbb{C}^{*})^{m}$,
the volume form on $M_P$ is given by
\[
\omega_{P}^{m}/m!=\frac{1}{(2\pi)^{m}}\det A(\rho) \,d\rho
d\theta\;.\] \label{volume}
\end{lem}

\begin{proof} We note that
\begin{equation}
A(\rho)=
\sum_{\beta \in P}k_{\beta}(\rho) \beta \otimes \beta -
\left(
\sum_{\beta \in P}k_{\beta}(\rho)\beta
\right)
\otimes
\left(
\sum_{\beta \in P}k_{\beta}(\rho)\beta
\right)={\rm
Hess}_{\rho} \log k(\rho),
\label{hess1}
\end{equation}
where $k_{\beta}(\rho)=
|\wh{m}_{\beta}^{P}(e^{\rho/2})|^{2}$ and $k(\rho)=\sum_{\be \in P}|c_{\be}|^{2}
e^{\ispa{\rho,\be}}$.
The conclusion follows from (\ref{hess1}),
recalling that
\begin{equation}
\label{volume0}
\omega_{P}=\Phi_{P}^{c*}\omega_{{\rm FS}} =
\frac{\sqrt{-1}}{2\pi}
\partial \bar\partial \log
\sum_{\beta}|c_{\beta}|^{2}|\chi_{\beta}(z)|^{2}=
\frac{\sqrt{-1}}{2\pi}
\partial \bar\partial \log
k(\rho).
\end{equation}
\end{proof}

It should be noted that $b_{x}^{P}(\rho)$ grows as $|\rho| \to \infty$, as stated in the following simple
lemma.

\begin{lem}
Let $K \subset P^{o}$ be a compact set. Then there exists positive constants $R>0$, $c>0$
such that $f(x,\rho) \geq c$ for $(x,\rho) \in K \times \R^{m}$, $|\rho| \geq R$.
\label{Lexpfunc2}
\end{lem}

\begin{proof}
For any $(x,\rho) \in P^{o} \times \R^{m}$, we define
\[
M(x,\rho)=\max_{\beta \in P \cap \Z^{m}}\ispa{\rho,\beta -x}.
\]
If $x \in P^{o}$, then the polytope $P -x$ contains the origin in its interior.
Thus, clearly we have $M(x,\rho) >0$ for any $(x,\rho) \in P^{o} \times (\R^{m}\setminus 0)$.
Next, we note that the function $(x,\rho) \mapsto M(x,\rho)$ is continuous.
To see this, let $(x_{n},\rho_{n})$ be a sequence such that $(x_{n},\rho_{n}) \to (x,\rho) \in P^{o} \times \R^{m}$.
Then, for any $\beta \in P \cap \Z^{m}$,
\[
|\ispa{\rho_{n},\beta-x_{n}}-\ispa{\rho,\beta-x}|
\leq C|\rho_{n}-\rho| +|\rho||x_{n}-x|.
\]
By using this inequality, we can show that
\[
|M(x_{n},\rho_{n})-M(x,\rho)| \leq  C|\rho_{n}-\rho| +|\rho||x_{n}-x|.
\]
Now, for a compact set $K \subset P^{o}$, we set
\[
M(K)=\min_{(x,\rho) \in K \times \R^{m},|\rho|=1}
M(x,\rho) >0.
\]
We set $c_{0}=\min_{\beta \in P\cap\Z^{m}}|c_{\beta}|^{2}$. Since $P\cap\Z^{m}$ is a finite set,
there exists $\beta=\beta(x,\rho)$ such that $M(x,\rho)=\ispa{\rho,\beta(x,\rho)-x}$ for $(x,\rho)$ with $x \in P^{o}$
and $|\rho|=1$.
Thus, for $x \in K \subset P^{o}$ and $\rho \neq 0$, we have
\[
e^{f(x,\rho)}=\sum_{\beta}|c_{\beta}|^{2}e^{\ispa{\rho,\beta-x}}\geq
c_{0}e^{|\rho|\ispa{\frac{\rho}{|\rho|},\beta(x,\frac{\rho}{|\rho|})-x}}
=c_{0}e^{|\rho|M(x,\frac{\rho}{|\rho|})}
\geq c_{0}e^{|\rho|M(K)}
\]
for $x \in K$ and $\rho \neq 0$, which completes that proof.
\end{proof}

\vspace{5pt}

\noindent{\it Completion of the proof of Theorem \ref{pointww}:\/}
Recalling that $z=e^{\rho/2 +i\theta}$, we write
$|\varphi_{\alpha_{N}}^{P}(\rho)|_{P}$ instead of
$|\varphi_{\alpha_{N}}^{P}(z)|_{P}$. By the definition of the
Hermitian metric on $L_{P}^{N}$, we have
\[
|\varphi_{\alpha_{N}}^{P}(\rho)|^{2}_{P}
=\frac{|\chi_{\alpha_{N}}^{P}(\rho)|_{P}^{2}}{\|\chi_{\alpha_{N}}^{P}\|^{2}}
=\frac{1}{\|\chi_{\alpha_{N}}^{P}\|^{2}}
\frac{e^{\ispa{\rho,\alpha_{N}}}}{|\Phi_{P}^{c}(e^{\rho/2})|^{2N}},
\]
where, in the right hand side, $|\Phi_{P}^{c}(e^{\rho/2})|$ denotes the usual norm in $\C^{d+1}$.
By the definition \eqref{PHPC} of the monomial embedding $\Phi_{P}^{c}$, we have
\[
|\Phi_{P}^{c}(e^{\rho/2})|^{2}=k(\rho),\quad \rho \in \R^{m}.
\]
Hence, we have
\[
|\varphi_{\alpha_{N}}^{P}(\rho)|_{P}^{2}=
\frac{e^{-Nf(\alpha_{N}/N,\rho)}}{\|\chi_{\alpha_{N}}^{P}\|^{2}}
=\frac{e^{-Nf(x,\rho)}R_{N}(x,\rho)}{\|\chi_{\alpha_{N}}^{P}\|^{2}},
\quad R_{N}(x,\rho)=e^{N\ispa{\rho,\alpha_{N}/N-x}},
\]
where the function $f(x,\rho)$ for $x \in P^{o}$ is defined in \eqref{binv}.
We note that, since $\alpha_{N}/N-x=O(N^{-1})$, we have
\begin{equation}
\label{high}
|\partial_{\rho}^{L}R_{N}(x,\rho)| \leq C_{L}R_{N}(x,\rho)
\end{equation}
for every multi-index $L$. By Lemma \ref{volume}, the
$\lcal^{2}$-norm of the un-normalized monomial
$\chi_{\alpha_{N}}^{P}$ is given by
\[
\|\chi_{\alpha_{N}}^{P}\|^{2}=
\int_{\R^{m}}e^{-Nf(\alpha_{N}/N,\rho)}\det A(\rho)\,d\rho.
\]
Here it should be noted that $\det A(\rho)$ is a positive integrable function on $\R^{m}$.
By Lemma \ref{Lexpfunc2}, we can choose $R>0$, $c>0$ such that $|\rho_{x}^{P}| <R$ and
$f(\alpha_{N}/N,\rho) \geq c$ for any $|\rho| \geq R$ and $N$.
Thus, by choosing a cut-off function $g(\rho)$ suitably, we may write
\[
\|\chi_{\alpha_{N}}^{P}\|^{2}
=e^{-Nf(x,\rho_{x}^{P})}
\int e^{-Nb_{x}^{P}(\rho)}R_{N}(x,\rho)g(\rho)\det A(\rho)\,d\rho +O(e^{-cN}).
\]
Recall that $b_{x}^{P}(\rho)=0$ if and only if $\rho
=\rho_{x}^{P}$, and that $\rho=\rho_{x}^{P}$ is the unique
critical point of $b_{x}^{P}$. The Hessian of $b_{x}^{P}$ at
$\rho=\rho_{x}^{P}$ is the positive definite symmetric matrix
$A(P,x)$. Thus, by the Morse lemma, there exists a change of
coordinates $\kappa$ from a neighborhood of the origin to a
neighborhood of $\rho_{x}^{P}$ such that $\kappa(0)=\rho_{x}^{P}$
and that
\[
b_{x}^{P} \comp \kappa(\xi) = \ispa{A(P,x)\xi,\xi}/2,\quad |\det D\kappa(0)|=1.
\]
By choosing the cut-off function $g$ suitably, we get
\[
\int e^{-Nb_{x}^{P}(\rho)}g(\rho)R_{N}(x,\rho)\det A(\rho)\,d\rho
=\int e^{-N\ispa{A(P,x)\xi,\xi}}G_{N}(x,\xi)\,d\xi,
\]
where $G_{N}(x,\xi)$ is a compactly supported function in $\xi$
such that $G_{N}(x,0)=\det A(P,x)$.
By \eqref{high}, the derivatives of $G_{N}(x,\xi)$ with respect to $\xi$
are all bounded uniformly in $N$.
Therefore, by using the Plancherel formula
and a formula for the Fourier transform of the Gaussian functions, we obtain
\[
\int e^{-Nb_{x}^{P}(\rho)}g(\rho)R_{N}(x,\rho)\det A(\rho)\,d\rho
=\left(
\frac{N}{2\pi}
\right)^{-m/2}
\sqrt{\det A(P,x)}(1+O(N^{-1})),
\]
which completes the proof.
\hfill\qedsymbol

\begin{rem} Since $\chi_{\alpha_{N}}^{P}$ is a  monomial, the
 asymptotics of $\|\chi_{\alpha_{N}}^{P}\|$ is essentially the
 same calculation as the asymptotics of the  $L^{2k}$ norm of another monomial.
 Thus, determining the  pointwise asymptotics of monomials is  equivalent to
 determining the asymptotics of their $L^{2k}$ norms.

\end{rem}

\subsection{Pointwise asymptotics: boundary lattice points}

Next, we consider the ray $\N\alpha$ for a lattice point $\alpha
\in P \cap\Z^{m}$, which is allowed to lie in the boundary
$\partial P$. In such a case, we need to work with other
coordinates than the usual coordinates on the open orbit
$(\C^{*})^{m}$ (\cite{SZ2}), since the open orbit $(\C^{*})^{m}$
does not cover the set $\mu_{P}^{-1}(\partial P)$.

In the following, we mean that the faces are disjoint, and the
facet is a face of codimension one. Thus we call the closed face
(or facet) $\bar{F}$ the closure of $F$ in the minimal affine
subspace containing $F$. To describe the coordinates, let $v_{0}$
be a vertex of $P$. Since our polytope $P$ is Delzant, we can
choose lattice points $\alpha^{1},\ldots,\alpha^{m}$ in $P$ such
that each $\alpha^{j}$ is in an edge incident to the vertex
$v_{0}$, and the vectors $v^{j}:=\alpha^{j}-v_{0}$ form a basis of
$\Z^{m}$. We choose (open) facets $F_{j}$, $j=1,\ldots,m$ incident
at $v_{0}$ so that $\alpha^{j} \not \in F_{j}$.

\begin{lem}
\label{divisor}
Let $\alpha \in P \cap \Z^{m}$, and $z \in M_{P}$. Then, $\chi_{\alpha}^{P}(z)=0$
if and only if
\begin{equation}
\label{Nulset}
\mu_{P}(z) \in \bigcup \{\bar{F}\,;\,F \mbox{ is a facet }\alpha \not \in \bar{F}\}.
\end{equation}
\end{lem}

\begin{proof}
If $z \in (\C^{*})^{m}$ then automatically we have $\chi_{\alpha}^{P}(z)\neq 0$. Thus,
we may assume that $z \in \mu_{P}^{-1}(\partial P)$.
First, assume that $\mu_{P}(z) \in \bar{F}$ for some closed facet $\bar{F}$ which does not contain $\alpha$.
The last formula in \eqref{muP} for the moment map is globally defined,
for the function $|\wh{m}_{\beta}^{P}|^{2}$ is globally defined. Since $\mu_{P}(z) \in \bar{F}$,
the coefficients in $\mu_{P}(z)$ of the lattice points $\beta \not \in \bar{F}$ must vanish.
Thus, we have $\chi_{\alpha}^{P}(z)=0$. Conversely, assume that $\mu_{P}(z) \in \partial P$
is not in the set described in \eqref{Nulset}.
In our convention, the faces are disjoint and the boundary $\partial P$ is the disjoint union
of faces. Thus, that $\mu_{P}(z) \in \partial P$ is not in the set described in \eqref{Nulset} is equivalent to say
that there exists an open face $E$ such that
$\alpha \in \bar{E}$ and $\mu_{P}(z) \in E$.
Let $v_{1},\ldots,v_{l}$ be the set of vertex of $\bar{E}$.
Since $v_{j}$ and $\alpha$ are lattice points, there exists a positive integer $n_{0}$ such that
$n_{0}\alpha=\sum_{j=1}^{l} n_{j}v_{j}$ with $n_{j}$ integer such that $\sum n_{j}=n_{0}$.
Thus we have
\[
(\chi_{\alpha}^{P})^{\otimes n_{0}}(z)=(\chi_{v_{1}}^{P})^{\otimes n_{1}}(z) \otimes
\cdots \otimes
(\chi_{v_{l}}^{P})^{\otimes n_{l}}(z).
\]
Since $\mu_{P}(z)$ is in the interior $E$ of the face (polytope) $\bar{E}$, each $\chi_{v_{j}}^{P}(z)$
can not vanish, and hence $\chi_{\alpha}^{P}(z) \neq 0$.
\end{proof}

Using Lemma \ref{divisor}, we set
\[
U_{v_{0}}:=\{z \in M_{P}\,;\,\chi_{v_{0}}^{P}(z) \neq 0\},
\]
which covers $M_{P}$ as $v_{0}$ varies over all vertices.
We define
\begin{equation}
\label{CChange}
\eta:(\C^{*})^{m} \to (\C^{*})^{m}, \quad
\eta(z)=\eta_{j}(z):=(z^{v^{1}},\ldots,z^{v^{m}}).
\end{equation}
The map $\eta$ is a diffeomorphism and the inverse is given by
\[
z:(\C^{*})^{m} \to (\C^{*})^{m},\quad z(\eta)=(\eta^{\Gamma
e^{1}},\ldots,\eta^{\Gamma e^{m}}),
\]
where $e^{j}$ is the standard basis for $\R^{m}$ (or for $\C^{m}$
over $\C$), and $\Gamma$ is an $m \times m$-matrix with $\det
\Gamma =\pm 1$ and integer coefficients defined by
\[
\Gamma v^{j}=e^{j},\quad v^{j}=\alpha^{j}-v_{0}.
\]
By definition, we have the obvious formula:
\[
\chi_{\alpha^{j}}^{P}(z)=\eta_{j}(z)\chi_{v^{0}}^{P}(z),\quad
z \in (\C^{*})^{m}.
\]
By Lemma \ref{divisor}, $\eta_{j}(z) \to 0$ if $z \in U_{v_{0}}$, $z \to \mu_{P}^{-1}(\bar{F}_{j})$.
Since $\alpha^{j} \not \in \bar{F}_{j}$, we have
\[
(\chi_{\alpha^{j}}^{P})^{-1}(0) \cap U_{v_{0}}
=\mu_{P}^{-1}(\bar{F}_{j}).
\]
The set $U_{v_{0}} \setminus (\C^{*})^{m}$ is the union of
the sets $\mu_{P}^{-1}(\bar{F}_{j})$, and hence the map $\eta$ extends a homeomorphism:
\[
\eta:U_{v_{0}} \to \C^{m},\quad \eta(z_{0})=0,\quad
z_{0}=\mbox{ the fixed point corresponding to } v_{0}.
\]
By this homeomorphism, the set $\mu_{P}^{-1}(\bar{F}_{j})$
corresponds to the set $\{\eta \in \C^{m}\,;\,\eta_{j}=0\}$.
This coordinate $\eta=(\eta_{1},\ldots,\eta_{m})$ is useful to explain toric subvarieties
corresponding to faces. Namely, let $\bar{F}$ be a closed face with $\dim F=m-r$ which contains $v_{0}$.
Since $v_{0} \in \bar{F}$, we can choose $F_{i_{1}},\ldots,F_{i_{r}}$ such that
$\bar{F}=\bar{F}_{i_{1}} \cap \cdots \cap \bar{F}_{i_{r}}$.
Then the subvariety $\mu_{P}^{-1}(\bar{F})$ corresponding $\bar{F}$ is expressed, in the coordinate
neighborhood $U_{v_{0}}$, by
\[
\mu_{P}^{-1}(\bar{F}) \cap U_{v_{0}}=
\{\eta \in \C^{m}\,;\,\eta_{i_{j}}=0,\quad j=1,\ldots,r\}.
\]

Now, we fix a lattice point $\alpha$ in a (relatively open)
face $F$ of dimension $\dim F=m-r$ such that $v_{0} \in \bar{F}$.
Without loss of generality, we may assume that $\bar{F}=\bar{F}_{1}\cap \cdots \cap \bar{F}_{r}$.

To state a result for the lattice point $\alpha$ in the boundary
corresponding to Theorem \ref{pointww}, we need to find a function
corresponding to the function $b_{\alpha}^{P}$. Since our
coordinate $\eta$ is based on the lattice points
$\alpha^{j}-v_{0}$, it is reasonable to introduce a new polytope
defined by the affine linear transformation
\[
\tilde{\Gamma}:\R^{m} \ni u \to \Gamma u -\Gamma v_{0} \in \R^{m},
\]
which maps $\Z^{m}$ bijectively onto itself.
We set $Q:=\tilde{\Gamma}(P)$. Then $Q$ is contained in the positive orthant $\{x \in \R^{m}\,;\,x_{j} \geq 0\}$,
and we have $\tilde{\Gamma}(\bar{F}_{j})=\{x \in Q\,;\,x_{j}=0\}$.
The face of $Q$ corresponding to $F$ is then given by
\[
Q_{F}:=\tilde{\Gamma}(\bar{F})=\{x \in Q\,;\, x_{j}=0,\quad j=1,\ldots,r\}.
\]
We denote a point in $U_{v_{0}} \cong \C^{m}$ as $\eta=(\xi,\zeta) \in \C^{m}=\C^{r} \times \C^{m-r}$.
In this expression, $\zeta =(0,\zeta)$ is a local coordinate of the submanifold $\mu_{P}^{-1}(\bar{F})$.
The modulus square of the monomials $|\chi_{N\alpha}^{P}|_{P}^{2}$ with $\alpha \in P \cap \Z^{m}$ is,
in this coordinate, given by
\begin{equation}
\label{monch}
\begin{gathered}
|\chi_{N\alpha}^{P}(z)|^{2}=
\frac{|\eta^{\tilde{\Gamma}(\alpha)}|^{2N}}{K(\eta)^{N}},\\
K(\eta)=\sum_{\gamma \in Q \cap \Z^{m}}a_{\gamma}|\eta^{\gamma}|^{2},
\quad a_{\gamma}=|c_{\tilde{\Gamma}^{-1}(\gamma)}|^{2}.
\end{gathered}
\end{equation}
We then introduce the `moment map' corresponding to the face $F$ by:
\begin{equation}
\label{momentF}
\mu_{F}:\R^{m-r} \to Q_{F},\quad
\mu_{F}(\rho)=\sum_{(0,\nu) \in Q_{F}}
\frac{a_{\nu}e^{\ispa{\rho,\nu}}}{k_{F}(\rho)}(0,\nu),
\end{equation}
where $a_{\nu}=|c_{\tilde{\Gamma}^{-1}(0,\mu)}|^{2}$,
and the function $k_{F}(\rho)$ is given by
\begin{equation}
\label{charF}
k_{F}(\rho):=\sum_{(0,\nu) \in Q_{F}}a_{\nu}e^{\ispa{\rho,\nu}}.
\end{equation}

As mentioned above, the submanifold $\mu_{P}^{-1}(\bar{F})$ has the coordinate
$\zeta \mapsto (0,\zeta) \in \C^{r} \times \C^{m-r}$, and the torus ${\bf T}^{m-r}$ acts on it.
Thus, it is natural to use the coordinate
\[
\zeta=e^{\rho/2 +i\theta},\quad \rho,\theta \in \R^{m-r}.
\]
Then, we write $\eta =(\xi,\zeta)=(\xi,\rho)$ for $\zeta=e^{\rho/2}$.
Since we have assumed $\alpha \in F$, we may write
\[
\tilde{\Gamma}(\alpha)=(0,\tilde{\alpha}) \in Q_{F}.
\]
We also define
\begin{equation}
\label{Phf}
\begin{gathered}
s_{\alpha}(\xi,\rho):=\log K(\xi,\rho)-\ispa{\rho,\tilde{\alpha}}
=\log (k_{F}(\rho)+\ell_{F}(\xi,\rho))-\ispa{\rho,\tilde{\alpha}},\\
\ell_{F}(\xi,\rho)=\sum_{(\mu,\nu) \in Q,\mu \neq 0}
a_{(\mu,\nu)}|\xi^{\mu}|^{2}e^{\ispa{\rho,\nu}}.
\end{gathered}
\end{equation}

\begin{prop}
\label{LOCALSS}
In the coordinate $\eta=(\xi,\zeta)$ as above, we write $|\varphi_{N\alpha}^{P}(\xi,\rho)|_{P}$
for the modulus square of the monomial $|\varphi_{N\alpha}^{P}|_{P}^{2}$.
We also write $\tilde{\Gamma}(\alpha)=(0,\tilde{\alpha})$, which is in the interior of the polytope $Q_{F}$.
Then we have
\[
|\varphi_{N\alpha}^{P}(\xi,\rho)|_{P}^{2}
=(2\pi)^{r}\left(
\frac{N}{2\pi}
\right)^{(m+r)/2}
\frac{e^{-N\Psi_{\alpha}(\xi,\rho)}}
{\sqrt{\det A(F,\alpha)}}(1+O(N^{-1})),
\]
where the function $\Psi_{\alpha}(\xi,\rho)$ is given by
\begin{equation}
\label{Phf2}
\Psi_{\alpha}(\xi,\rho)=s_{\alpha}(\xi,\rho)-s_{\alpha}(0,\rho_{\alpha}^{F})\quad
\rho_{\alpha}^{F}=\mu_{F}^{-1}(\tilde{\alpha}) \in \R^{m-r}.
\end{equation}
and $(m-r) \times (m-r)$ positive definite matrix $A(F,\alpha)$ is given by
\begin{equation}
\label{matrixAF}
A(F,\alpha)=\sum_{(0,\nu) \in Q_{F}}
\frac{a_{\nu}e^{\ispa{\rho_{\alpha}^{F},\nu}}}{k_{F}(\rho_{\alpha}^{F})}
\nu \otimes \nu -\tilde{\alpha} \otimes \tilde{\alpha}.
\end{equation}
\end{prop}

We shall prove Proposition \ref{LOCALSS} in the rest of this subsection.
First of all, we need the following simple lemma:

\begin{lem}
\label{volume22}
In the coordinate $(\xi,\zeta=e^{\rho/2 +i\theta})$ on $\C^{r} \times (\C^{*})^{m-r} \subset U_{v_{0}}$,
the volume form $\omega_{P}^{m}/m!$ is given by
\begin{equation}
\label{DETV}
\frac{\omega_{P}^{m}}{m!}=
\frac{1}{\pi^{r}(2\pi)^{m-r}}
L(\xi,\rho)\,dm(\xi)d\rho d\theta,
\end{equation}
where $dm(\xi)$ denotes the Lebesgue measure on $\C^{r}$,
and the function $L(\xi,\rho)$ is given by the determinant of the following
$m \times m$ matrix:
\begin{equation}
\label{DETM}
L(\xi,\rho)=\det
\begin{pmatrix}
\frac{\partial^{2}\log K}{\partial \xi \partial \bar{\xi}} &
\frac{\partial^{2}\log K}{\partial \xi \partial \rho} \\
\frac{\partial^{2}\log K}{\partial \rho \partial \bar{\xi}} &
\frac{\partial^{2}\log K}{\partial \rho \partial \rho}
\end{pmatrix},
\end{equation}
where $K(\xi,\rho)$ is the function defined in \eqref{monch}.
\end{lem}

\begin{proof}
For the function $f(\xi,\zeta)$ on $\C^{r} \times (\C^{*})^{m-r}$ independent of the
variable $\theta$ in $\zeta=e^{\rho/2 +i\theta}$, then the derivatives $(\partial f)/(\partial \zeta_{j})$,
$(\partial f)/(\partial \bar{\zeta}_{j})$ is given, respectively, by
\[
\frac{\partial f}{\partial \zeta_{j}}=
\frac{1}{\zeta_{j}}\frac{\partial f}{\partial \rho_{j}},\quad
\frac{\partial f}{\partial \bar{\zeta}_{j}}=
\frac{1}{\bar{\zeta}_{j}}\frac{\partial f}{\partial \rho_{j}}.
\]
We write $\eta=(\xi,\zeta)$. Then, by this relation, we have
\[
\det
\left(
\frac{\partial^{2}\log K}{\partial \eta_{j}\partial \bar{\eta}_{k}}
\right)
=\left(
\prod_{j=1}^{m-r}|\zeta_{j}|^{2}
\right)^{-1}
L(\xi,\rho),
\]
where $L(\xi,\rho)$ is given by \eqref{DETM}. But, we have
\[
dm(\zeta)=\frac{1}{2^{m-r}}
\left(
\prod_{j=1}^{m-r}|\zeta_{j}|^{2}
\right)
d\rho d\theta,
\]
where $dm(\zeta)$ denotes the Lebesgue measure on $\C^{m-r}$.
Combining this with \eqref{volume0}, we obtain the assertion.
\end{proof}

The following lemma can be shown by the same argument as in the proof
of Lemma \ref{posmat}, and we shall omit the proof.

\begin{lem}
\label{posmatF}
The $(m-r) \times (m-r)$ matrix defined by
\[
A_{F}(\rho):=\partial \mu_{F}(\rho)=
\sum_{(0,\nu) \in Q_{F}}\frac{a_{\nu}e^{\ispa{\rho,\nu}}}{k_{F}(\rho)}
\nu \otimes \nu -\mu_{F}(\rho) \otimes \mu_{F}(\rho)
\]
is positive definite for every $\rho \in \R^{m-r}$.
\end{lem}

Note that, the map $\mu_{F}:\R^{m-r} \to Q_{F}^{o}$ is a diffeomorphism, and the lattice point $\tilde{\alpha}$
is in $Q_{F}^{o}$. Thus, the vector $\rho_{\alpha}^{F}=\mu_{F}^{-1}(\tilde{\alpha})$ is well-defined.
Hence, the $(m-r) \times (m-r)$ matrix
\[
A(F,\alpha)=A_{F}(\rho_{\alpha}^{F})
\]
is positive definite.

\vspace{10pt}

\noindent{\it Completion of proof of Proposition \ref{LOCALSS}.}\hspace{3pt}The modulus square of
the monomial $|\varphi_{N\alpha}^{P}|_{P}^{2}$ in this coordinate is given by
\[
\begin{split}
|\varphi_{N\alpha}^{P}(\xi,\rho)|_{P}^{2} &=
\frac{|\chi_{N\alpha}^{P}(\xi,\rho)|_{P}^{2}}{\|\chi_{N\alpha}^{P}\|^{2}} \\
& =\frac{1}{\|\chi_{N\alpha}^{P}\|^{2}}
e^{-Ns_{\alpha}(\xi,\rho)},
\end{split}
\]
where the function $s_{\alpha}(\xi,\rho)$ is defined by
\eqref{Phf}. Thus, as in the proof of Theorem \ref{pointww}, what
we need to analyze is the $\lcal^{2}$-norm
\begin{equation}
\label{L2norm1}
\|\chi_{N\alpha}^{P}\|^{2}
=\frac{1}{\pi^{r}}
\int_{\C^{r} \times \R^{m-r}}
e^{-Ns_{\alpha}(\xi,\rho)}
L(\xi,\rho)\,dm(\xi)d\rho,
\end{equation}
where we have used Lemma \ref{volume22}.
We note that the function $\ell_{F}(\xi,\rho)$ defined in \eqref{Phf} is of the form
\begin{equation}
\label{expand}
\begin{gathered}
\ell_{F}(\xi,\rho)=\sum_{k=1}^{r}f_{k}(\rho)|\xi_{k}|^{2}+r(\xi,\rho),\\
f_{k}(\rho)=\sum_{(e_{k}^{r},\nu) \in Q}
a_{(e_{k}^{r},\nu)}e^{\ispa{\rho,\nu}},\quad
r(\xi,\rho)=\sum_{(\mu,\nu) \in Q,|\mu| \geq 2}a_{(\mu,\nu)}|\xi^{\mu}|^{2}e^{\ispa{\rho,\nu}},
\end{gathered}
\end{equation}
where $e_{k}^{r}$ are the standard basis for $\R^{r}$.
The function $r(\xi,\rho)$ is of order $\geq 4$ in $\xi$, and hence its derivative
up to the second order vanish at $\xi=0$.
Thus, the function $s_{\alpha}(\xi,\rho)$ has only one critical point $(\xi,\rho)=(0,\rho_{\alpha}^{F})$.
It is not hard to show that the Hessian $Hs_{\alpha}(0,\rho_{\alpha})$ of the function \eqref{Phf}
at the critical point $(0,\rho_{\alpha}^{F})$ is given by
\begin{equation}
\label{hess222}
Hs_{\alpha}(0,\rho_{\alpha}^{F})=
\begin{pmatrix}
\frac{2f_{1}(\rho_{\alpha}^{F})}{k_{F}(\rho_{\alpha}^{F})} & & & & & & & \\
& \ddots & & & & &\\
& & \frac{2f_{r}(\rho_{\alpha}^{F})}{k_{F}(\rho_{\alpha}^{F})} & & & & & \\
& & & \frac{2f_{1}(\rho_{\alpha}^{F})}{k_{F}(\rho_{\alpha}^{F})} & & & \\
& & & & \ddots & & \\
& & & & & \frac{2f_{r}(\rho_{\alpha}^{F})}{k_{F}(\rho_{\alpha}^{F})} & \\
& & & & & & A(F,\alpha)
\end{pmatrix},
\end{equation}
where the $(m-r) \times (m-r)$-matrix $A(F,\alpha)$ is given by
\eqref{matrixAF}, and we have used the coordinate
$(x_{j},y_{j},\rho)$ with $\xi_{j}=x_{j}+iy_{j}$. Here, it should
be noted that the polytope $Q=\tilde{\Gamma}(P)$ contains the
standard basis in $\R^{m}=\R^{r} \times \R^{m-r}$. Thus the
lattice points $(e_{j}^{r},0)$ with the standard basis $e_{j}^{r}$
in $\R^{r}$ is in the polytope $Q$, and hence the functions
$f_{j}$ are all positive. This combined with Lemma \ref{posmatF}
shows that the Hessian $Hs_{\alpha}(0,\rho_{\alpha})$ is positive
definite. By the same argument as in the proof of Lemma
\ref{Lexpfunc2}, the function $k_{F}(\rho)+\ell_{F}(\xi,\rho)$
tends to $\infty$ as $|\xi|+|\rho| \to \infty$. Therefore, by the
standard Laplace method as in the proof of Theorem \ref{pointww},
we have
\begin{equation}
\label{laplace1}
\|\chi_{N\alpha}^{P}\|^{2}
=\frac{1}{\pi^{r}}
\left(
\frac{N}{2\pi}
\right)^{-(m+r)/2}
\frac{e^{-Ns_{\alpha}(0,\rho_{\alpha}^{F})}}
{\sqrt{\det Hs_{\alpha}(0,\rho_{\alpha}^{F})}}
L(0,\rho_{\alpha}^{F})(1+O(N^{-1})).
\end{equation}
A direct computation will show that
\begin{equation}
\label{DET}
\det Hs_{\alpha}(0,\rho_{\alpha}^{F})
=\frac{\det A(F,\alpha)}{k_{F}(\rho_{\alpha}^{F})^{2r}}
\left(
\prod_{j=1}^{r}2f_{j}(\rho_{\alpha}^{F})
\right)^{2},\quad
L(0,\rho)=
\frac{\det A_{F}(\rho)}{k_{F}(\rho_{\alpha}^{F})^{r}}
\left(
\prod_{j=1}^{r}f_{j}(\rho)
\right),
\end{equation}
and hence, the asymptotics \eqref{laplace1} can be written in the form:
\[
\|\chi_{N\alpha}^{P}\|^{2}
=\frac{1}{(2\pi)^{r}}
\left(
\frac{N}{2\pi}
\right)^{-(m+r)/2}
\sqrt{\det A(F,\alpha)}
e^{-Ns_{\alpha}(0,\rho_{\alpha}^{F})}(1+O(N^{-1})).
\]
Dividing $|\chi_{N\alpha}|^{2}_{P}$ by the above, we conclude the assertion.
\hfill\qedsymbol

When our fixed lattice point $\alpha$ is a vertex, say $\alpha=v_{0}$ in the description
of the coordinate $\eta$,
the matrix $A(F,\alpha)$ is not defined suitably.
However, clearly the similar asymptotics can be deduce by the same method.

\begin{prop}
\label{LOCALV}
Suppose that $\alpha$ is a vertex of the polytope $P$. Then we have
\[
|\varphi_{N\alpha}^{P}(\eta)|_{P}^{2}
=\left(
\frac{N}{|c_{\alpha}|^{2}}
\right)^{m}
e^{-N(\log K(\eta) -\log |c_{\alpha}|^{2})}(1+O(N^{-1})),
\]
where the function $K(\eta)$ on $\C^{m}$ is given by \eqref{monch}. We also have
$|c_{\alpha}|^{2}=|\chi_{\alpha}^{P}(z_{\alpha})|^{-2}$ where $z_{\alpha}$ is the
fixed point for the Hamiltonian ${\bf T}^{m}$-action such that $\mu_{P}(z_{\alpha})=\alpha$.
\end{prop}

\begin{proof}
In the description of the coordinate $\eta$, we put $\alpha =v_{0}$.
Then, clearly we have $K(0)=|c_{\alpha}|^{2}=|\chi_{\alpha}^{P}(z_{\alpha})|^{-2}$, where
the fixed point $z_{\alpha}$ corresponds to the origin $\eta=0$ and
the function $K(\eta)$ is defined in \eqref{monch}.
In this case, we just use the coordinate $\eta$ itself without change of variable.
In this coordinate, the monomial $|\chi_{N\alpha}^{P}|_{P}^{2}$ is given by
\[
|\chi_{N\alpha}^{P}(\eta)|_{P}^{2}
=e^{-N\log K(\eta)}.
\]
The volume measure $\omega_{P}^{m}/m!$ is of the form:
\[
\frac{\omega_{P}^{m}}{m!}
=\frac{1}{\pi^{m}}
\det L(\eta)dm (\eta),\quad
L(\eta)=
\left(
\frac{\partial^{2}\log K(\eta)}{\partial \eta_{j} \partial \bar{\eta}_{k}}
\right)dm(\eta)
\]
with the Lebesgue measure $dm(\eta)$ on $\C^{m}$.
In is straight forward to see that the critical point of the function $\log K(\eta)$ is the origin,
and the determinant of the Hessian at the origin is given by
\[
H(\log K)(0)=
\begin{pmatrix}
\frac{2L(0)}{K(0)} & 0 \\
0 & \frac{2L(0)}{K(0)}
\end{pmatrix}.
\]
By using these facts with the Laplace method, we obtain the assertion.
\end{proof}

\subsection{Moments and $\lcal^{2k}$ norms}

As an application of the pointwise estimates,  one can prove that
eigenfunctions  `localize on tori'. We also determine $\lcal^{2k}$
norm of the monomials. The following is easily shown by using
Theorem \ref{pointww},  Propositions \ref{LOCALSS} and
\ref{LOCALV}, and the argument is the same as in their proofs, and
hence we shall omit the proof (see also the proof of Theorem
\ref{ml2knormP}).

\begin{prop}\label{lT}
\begin{itemize}

\item[{\rm (i)}] Let $\alpha_{N} \in NP \cap \Z$ be a sequence of lattice points
satisfying \eqref{GENL} for some point $x \in P^{o}$.
Then, the measure $|\varphi_{\alpha_{N}}^{P}|_{P}^{2}d\vol_{M_{P}}$ weak$^{*}$-converges
to the normalized Haar measure on the $m$-dimensional torus $\mu_{P}^{-1}(x)$, i.e.,
\[
\int_{M_{P}}\sigma |\varphi_{\alpha_{N}}^{P}|_{P}^{2}d\vol_{M_{P}}
\to \int_{\mu_{P}^{-1}(x)}\sigma d\theta
\]
for $\sigma \in \ccal (M_{P})$.

\item[{\rm (ii)}] Let $\alpha$ be a lattice point in $P$ with $\dim \mu_{P}^{-1}(\alpha)=m-r$.
Then, we have
\[
{\rm w}^{*}\!\mbox{-}\!\lim_{N \to \infty}|\varphi_{N\alpha}^{P}|_{P}^{2}
=d\theta_{\mu_{P}^{-1}(\alpha)},
\]
where $d\theta_{\mu_{P}^{-1}(\alpha)}$ denotes the normalized Haar measure
on the $(m-r)$-dimensional torus $\mu_{P}^{-1}(\alpha) \cong {\bf T}^{m-r}$.
\end{itemize}
\end{prop}

In the above proposition, if $\alpha$ is a vertex, then the left hand side denotes the Dirac measure
at the fixed point $z_{\alpha} \in M_{P}$ of the Hamiltonian ${\bf T}^{m}$-action corresponding to $\alpha$.

We note that $|\wh\phi_{N \alpha}^{P}|$ is invariant under the ${\bf T}^m$
action. We denote the Hilbert space of ${\bf T}^m$ invariant
functions by $\lcal^2_{{\rm inv}}(M_P)$. We can restate the conclusion
as follows: if  $\sigma \in C^{\infty}_{{\rm inv}}(M_{P})$, then
we can regard it as a function on the polytope $P$. We can also
regard $|\wh\phi_{N \alpha}^{P}|^2 $ as a function, say $|\wt\phi_{\alpha}^P(I)|^2$,
on $P$, equipped with action variables $I$. We then have:

\begin{cor} For any $\alpha \in P \cap \Z^{m}$, we have $|\wt\phi^P_{N
\alpha}|^2 dI \to \delta_{\alpha}$
; i.e., $\int_{P}
\sigma |\wt\phi_{\alpha}^P|^2\,  dI  \to \sigma(\alpha)$, for  $\sigma \in \ccal(P)$.
For $\alpha_{N} \in NP \cap\Z^{m}$ satisfying \eqref{GENL} with a point $x \in P^{o}$, we have
$|\wt{\varphi}_{\alpha_{N}}^{P}|^{2}dI \to \delta_{x}$.
\end{cor}

Next, we determine the asymptotics of the $\lcal^{2k}$-norm of the $\lcal^{2}$-normalized monomials.

\begin{theo}\label{ml2knormP}

\begin{itemize}

\item[{\rm (i)}] Let $\alpha_{N} \in NP \cap\Z^{m}$ be a sequence of lattice points satisfying \eqref{GENL}
with a point $x \in P^{o}$. Let $\|\varphi_{\gamma}^{P}\|_{2k}$ denote the $\lcal^{2k}$-norm
of the $\lcal^{2}$-normalized monomial $\varphi_{\gamma}^{P}$ with the weight $\gamma \in NP\cap\Z^{m}$.
Then we have
\begin{equation}
\label{Nonlat}
\|\varphi_{\alpha_{N}}^{P}\|_{2k}^{2k}=
\frac{1}{k^{m/2}}
\left(
\frac{N}{2\pi}
\right)^{m/2}
\frac{1}{(\det A(P,x))^{(k-1)/2}}(1+O_{k}(N^{-1})),
\end{equation}
where the $O_{k}(N^{-1})$ depends on $k$.

\item[{\rm (ii)}] Let $\alpha \in P$ be a lattice point with
$\dim \mu_{P}^{-1}(\alpha)=m-r$.
Then, for $r \leq m-1$, we have
\begin{equation}
\label{NonV}
\|\wh\phi^P_{N\alpha}\|_{2k}^{2k}
=
\frac{1}{k^{(m+r)/2}}
\left(
\frac{N}{2\pi}
\right)^{(k-1)(m+r)/2}
\left(
\frac{(2\pi)^{r}}{\sqrt{\det A(F,\alpha)}}
\right)^{k-1}
(1+O_{k}(N^{-1})).
\end{equation}
For a vertex $\alpha$, we have
\begin{equation}
\label{Vertex}
\|\wh\phi^P_{N\alpha}\|_{2k}^{2k}
=\frac{1}{k^{m}}
\left(
\frac{N}{|c_{\alpha}|^{2}}
\right)^{(k-1)m}(1+O_{k}(N^{-1})).
\end{equation}
\end{itemize}
\end{theo}

\begin{proof}
The proof of \eqref{Nonlat} and \eqref{Vertex} is the same as that
for \eqref{NonV}, so we shall give a proof of \eqref{NonV} only.
By Theorem \ref{pointww} and Proposition \ref{LOCALSS}, we have
\[
\begin{split}
\|&\wh\phi^P_{N\alpha}\|_{2k}^{2k}\\
&=
\frac{(2\pi)^{kr}}{\pi^{r}}
\left(
\frac{N}{2\pi}
\right)^{k(m+r)/2}\!\!\!
\frac{1}{(\det A(F,\alpha))^{k/2}}
\int_{\C^{r} \times \R^{m-r}}\!\!\!
e^{-Nk\Psi_{\alpha}(\xi,\rho)}L(\xi,\rho)\,dm(\xi)d\rho
(1+O_{k}(N^{-1})),
\end{split}
\]
where $\Psi_{\alpha}$ is given by \eqref{Phf2}.
As in the proof of Proposition \ref{LOCALSS}, the critical point of $\Psi_{\alpha}$
is only the point $(0,\rho_{\alpha}^{F})$ with $\rho_{\alpha}^{F}=\mu_{F}^{-1}(\tilde{\alpha})$,
$\tilde{\Gamma}(\alpha)=(0,\tilde{\alpha})$, and which is non-degenerate.
We have $\Psi_{\alpha}(0,\rho_{\alpha}^{F})=0$.
Thus, by the standard Laplace method, we have
\[
\|\wh\phi^P_{N\alpha}\|_{2k}^{2k}
=\frac{(2\pi)^{kr}}{\pi^{r}}
\left(
\frac{N}{2\pi}
\right)^{(k-1)(m+r)/2}
\frac{1}{(\det A(F,\alpha))^{k/2}}
\frac{L(0,\rho_{\alpha})}{\det H\Psi_{\alpha}(0,\rho_{\alpha})}(1+O_{k}(N^{-1})).
\]
Therefore, the assertion follows from \eqref{DET} in the proof of Proposition \ref{LOCALSS}.
\end{proof}

By using Proposition \ref{LOCALSS}, we can determine the limit of the sup-norm
$\|\varphi_{N\alpha}^{P}\|_{\infty}^{2}$.

\begin{prop}
\label{sup}
\begin{itemize}

\item[{\rm (i)}] Let $\alpha_{N} \in NP\cap\Z^{m}$ satisfy \eqref{GENL} with a point $x \in P^{o}$.
Then we have
\begin{equation}
\label{limN}
\lim_{N\to \infty}
\left(
\frac{N}{2\pi}
\right)^{m/2}\|\varphi_{\alpha_{N}}^{P}\|_{\infty}^{2}
=\frac{1}{\sqrt{\det A(P,x)}}.
\end{equation}

\item[{\rm (ii)}] Let $\alpha$ be a lattice point with $\dim \mu_{P}^{-1}(\alpha)=m-r$, $r \leq m-1$.
Let $(\xi_{N},\rho_{N})$ be the point where $|\varphi_{N\alpha}^{P}|_{P}^{2}$ attains its maximum.
Then, we have
\begin{equation}
\label{limitS}
\lim_{N \to \infty}
\left(
\frac{N}{2\pi}
\right)^{-(m+r)/2}
\|\varphi_{N\alpha}^{P}\|_{\infty}^{2}
=\frac{(2\pi)^{r}}{\sqrt{\det A(F,\alpha)}}.
\end{equation}
For a vertex $\alpha \in P$, we have
\begin{equation}
\label{limitV}
\lim_{N\to \infty}N^{-m}\|\varphi_{N\alpha}^{P}\|_{\infty}^{2}
=|c_{\alpha}|^{-2m}=
|\chi_{\alpha}^{P}(z_{\alpha})|_{P}^{2m},
\end{equation}
where $z_{\alpha} \in M_{P}$ is the unique fixed point
for the ${\bf T}^{m}$ such that $\mu_{P}(z_{\alpha})=\alpha$.
\end{itemize}
\end{prop}

\begin{proof}
We only give a proof of \eqref{limitS}.
The function $\Psi_{\alpha}(\xi,\rho)$ defined in \eqref{Phf2}
attains its minimum at the point $(0,\rho_{\alpha}^{F})$, and it tends to $\infty$ as $|\xi|+|\rho| \to \infty$.
Thus, the function $(\frac{N}{2\pi})^{-(m+r)/2}|\varphi_{N\alpha}^{P}|_{P}^{2}$ is of order $O(e^{-cN})$ outside
a compact neighborhood of $(0,\rho_{\alpha}^{F})$.
On a compact neighborhood $B$ of $(0,\rho_{\alpha}^{F})$, we have
\[
|\varphi_{N\alpha}^{P}(0,\rho_{\alpha}^{F})|_{P}^{2}
\leq \sup_{(\xi,\rho) \in B}|\varphi_{N\alpha}^{P}(\xi,\rho)|_{P}^{2}
\leq
\left(
\frac{N}{2\pi}
\right)^{(m+r)/2}
\frac{(2\pi)^{r}}{\sqrt{\det A(F,\alpha)}}
(1+O(N^{-1})).
\]
Now \eqref{limitS} follows from the above inequality.
The same argument with Proposition \ref{LOCALSS} shows \eqref{limitV}.
\end{proof}

\subsection{Asymptotics on projective space}
\label{projective}

The values of the $\lcal^{2k}$ norm in Proposition \ref{exact}
and in the projective-space case of Theorem \ref{ml2knorm} may seem to be different.
However, we can check that these two coincide by noting the following simple lemma.

\begin{lem}
The determinant $\det A(p\Sigma,\alpha)$ of the matrix $A(p\Sigma,\alpha)$ is
given by
\[
\det A(p\Sigma,\alpha)
=\frac{(p-|\alpha|)\alpha_{1}\cdots \alpha_{m}}{p}.
\]
\label{det}
\end{lem}

\begin{proof}
By applying the differential operator $x_{j}\partial_{x_{j}}$ twice
to the formula $(1+\sum_{l=1}^{m}x_{l})^{p}=\sum_{|\beta|\leq p}{p\choose\beta} x^{\beta}$
for a vector $x \in \mathbb{R}^{m}$, we have
\[
p(p-1)x_{i}x_{j}(1+\sum_{l}x_{l})^{p-2}+px_{j}(1+\sum_{l}x_{l})^{p-1}\delta_{ij}
=\sum_{|\beta| \leq p}{p\choose\beta} x^{\beta}\beta_{i}\beta_{j},
\]
where $\delta_{ij}$ is the Kronecker's delta, and we have set
${p\choose\beta}=p!/\beta !(p-|\beta|)!$. We put $r=\sum
e^{(\rho_{\alpha})_{j}}$. Then, by \eqref{special}, we have
$e^{\rho_{\alpha}}=\frac{1+r}{p}\alpha$ and
$|\widehat{m}^{p\Sigma}_{\beta}(e^{\rho_{\alpha}/2})|^{2}=e^{\ispa{\rho_{\alpha},\beta}}{p\choose\beta}
/(1+r)^{p}$. From this the second equation follows.
Substituting $\frac{1+r}{p}\alpha$ for $x$ in the
above formula, we obtain
\[
A(p\Sigma,\alpha)_{i,j}=\sum_{|\beta| \leq p} {p\choose\beta}
\frac{e^{\ispa{\rho,\beta}}}{(1+r)^{p}}\beta_{i}\beta_{j}-\alpha_{i}\alpha_{j}
= \alpha_{j}\delta_{ij} -\frac{1}{p}\alpha_{i}\alpha_{j}.
\]
We set $D_{m}(\alpha_{1},\ldots,\alpha_{m})=\det A(p,\alpha)$ with
$\alpha=(\alpha_{1},\ldots,\alpha_{m})$. Then a simple computation
shows that
\[
\frac{D_{m}(\alpha_{1},\ldots,\alpha_{m})}{\alpha_{1}\cdots\alpha_{m}}
=\frac{D_{m-1}(\alpha_{2},\ldots,\alpha_{m})}{\alpha_{2}\cdots\alpha_{m}}-\frac{1}{p}\alpha_{1}.
\]
Thus the lemma follows by  induction on the dimension $m$.
\end{proof}

Note that, for the simplex $p\Sigma$, the faces containing the origin is again a simplex in lower dimensional
vector space. Therefore, Lemma \ref{det} can be applied to compute $\det A(F,\alpha)$.

In the case of projective monomials, the
function $b_{\alpha}^{p\Si}$ can be expressed as follows.

\begin{prop}
For $z \in (\mathbb{C}^{*})^{m}$,
\[
b_{\alpha}^{p\Si}(z)=p\log(1+|z|^{2})-\log|z^{\alpha}|^2 +
\ispa{\widehat{\alpha},\log \widehat{\alpha}}-p\log p,
\]
where $\log \widehat{\alpha}=(\log
\widehat{\alpha}_{0},\ldots,\log \widehat{\alpha}_{m})$.
\end{prop}

This is a direct consequence of the formulas
$e^{\rho_{\alpha}}=\frac{1+r}{p}\alpha$, $r=\sum
e^{(\rho_{\alpha})_{j}}$ and $(1+\sum_{l}x_{l})^{p}=\sum_{|\beta|
\leq p}{p\choose\beta}x^{\beta}$ as mentioned before.

The pointwise asymptotics of the $\lcal^{2}$-normalized projective monomials is given in
the following proposition.

\begin{prop}
\label{LOCALP}
\begin{enumerate}

\item[{\rm (1)}] For $z \in (\mathbb{C}^{*})^{m}$,
\[
|\wh\phi_{N\alpha}^{p\Si}(z)|^{2} = \left( \frac{N}{2\pi} \right)^{m/2}
\frac{p^{1/2}e^{-Nb_{\alpha}^{p\Si}(z)}}{\sqrt{(p-|\alpha|)\alpha_{1}
\cdots \alpha_{m}}}(1+O(N^{-1})),
\]
uniformly $(\mathbb{C}^{*})^{m}$.

\item[{\rm (2)}] For  $z =e^{(\rho_{\alpha}+u/\sqrt{N})/2 +i\theta}$,
we have
\[
|\wh\phi_{N\alpha}^{p\Si}(z)|^{2} = \left( \frac{N}{2\pi} \right)^{m/2}
\frac{p^{1/2}e^{-(\ispa{\Delta(\alpha)u,u}-\frac{1}{p}\ispa{\alpha,u}^{2})/2}}
{\sqrt{(p-|\alpha|)\alpha_{1} \cdots \alpha_{m}}}(1+O(N^{-1/2}))
\]
uniformly for $|u| \leq c$, where $\Delta(\al)$ is the diagonal matrix
with entries $\alpha_{1},\ldots,\alpha_{m}$.
\end{enumerate}
\end{prop}

The assertion $(1)$ in the above is a restatement of Theorem
\ref{pointww} for the projective monomials. The assertion $(2)$
follows from $(1)$ and a Taylor expansion of the function
$b_{\alpha}^{p\Sigma}$.

\section{Asymptotics of distribution functions}
\label{invmprob}

In this section, we find asymptotics of rescaled and un-rescaled
distribution functions. Fix a lattice point $\alpha$ in $P$,
and  let $r$ denote the codimension of the face of $P$ (possibly
the open face $P^o$) containing $\al$. In analogy with
\eqref{CONSTX}, we define the constant $c(P,\alpha)$ by
\begin{equation}
\label{CONSTC}
c(P,\alpha):=
\left\{
\begin{aligned}
\frac{(2\pi)^{r}}{\sqrt{\det A(F,\alpha)}}\qquad \qquad
 &\quad \mbox{ if } r<m\\
\frac{(2\pi)^{m}}{|c_{\alpha}|^{2m}}=
(2\pi |\chi_{N\alpha}^{P}(z_{\alpha})|^{2}_{P})^{m}
&\quad \mbox{ if }\alpha \mbox{ is a vertex}.
\end{aligned}
\right.
\end{equation}

In the following discussion, we give the details for the case
where $\alpha_{N}=N\alpha$ with a lattice point $\alpha \in P
\cap\Z^{m}$. For general $\alpha_{N} \in NP\cap\Z^{m}$ satisfying
\eqref{GENL} for a point $x \in P^{o}$, one needs only to put
$r=0$ and replace $N\alpha$ and $\alpha$ by $\alpha_{N}$ and
$x$.

\subsection{Rescaled distribution functions}

We would like to understand the limit distributions of the measures $|\wh\phi^{P}_{N\alpha}|^{2}_{*}d\vol_{M_{P}}$,
namely, the limit of its distribution function
\begin{equation}
D_{N\alpha}(t):=\vol_{M_{P}}\{z\,;\,|\wh\phi^{P}_{N\alpha}(z)|^{2} >t\}.
\label{df}
\end{equation}
However, by Theorem \ref{ml2knorm}, the $k$-th moments of the
measure $|\wh\phi^{P}_{N\alpha}|^{2}_{*}d\vol_{M_{P}}$ tends to infinity as $N$ tends to infinity.
Therefore, we need to re-normalize the monomials.

We write
\begin{equation}
\label{renorm1}
{\rm dv}_{N}^{r}=
\left (
\frac{N}{2\pi}
\right )^{(m+r)/2}d\vol_{M_{P}}, \quad
f_{N\alpha}=
\left(
\frac{N}{2\pi}
\right)^{-(m+r)/4}
\varphi_{N\alpha}^{P}(z)
\end{equation}
so that
\[
\int_{M_{P}}
|f_{N\alpha}(z)|_{P}^{2}\,{\rm dv}_{N}^{r}(z)
=1.
\]
By Propositions \ref{limitS} and \ref{limitV},
we know that $\lim_{N \to \infty}\|f_{N\alpha}\|^{2}_{\infty}$ exists
and we have
\begin{equation}
\lim_{N \to \infty}\|f_{N\alpha}\|^{2}_{\infty}
=c(P,\alpha).
\label{const2}
\end{equation}
Furthermore, by Theorem \ref{ml2knorm}, we have
\begin{equation}
\|f_{N\alpha}\|_{\lcal^{2k}({\rm dv}^{r}_{N})}^{2k}=
\left(
\tf{N}{2\pi}
\right)^{-(m+r)(k-1)/2}\|\wh\phi^{P}_{N\alpha}\|_{2k}^{2k}
=\frac{c(P,\alpha)^{k-1}}{k^{(m+r)/2}}(1+O(N^{-1}))\;.
\label{renorm3}
\end{equation}
We consider the limit distribution of the sequence of  measures
\begin{equation}
\nu_{N,r}:=|f_{N\alpha}|^{2}_{*}{\rm dv}^{r}_{N}
\label{measure1}
\end{equation}
on the real line.
The distribution function $F_{N\alpha}^{r}(t)$ of the measure $\nu_{N,r}$ is given by
\begin{equation}
F_{N\alpha}^{r}(t)=
\left(
\tf{N}{2\pi}
\right)^{(m+r)/2}
\vol_{M_{P}} \{z \in M_{P}\,;\,|\wh\phi^{P}_{N\alpha}(z)|^{2} >
\left(
\tf{N}{2\pi}
\right)^{(m+r)/2}t
\}.
\label{df2}
\end{equation}
The distribution function $F_{N\alpha}^{r}$ defined above can be expressed, in terms of the
distribution function $D_{N\alpha}$ for the measure $|\wh\phi^{P}_{N\alpha}|^{2}_{*}d\vol_{M_{P}}$, as
\[
F^{r}_{N\alpha}(t)=
\left(
\tf{N}{2\pi}
\right)^{(m+r)/2}
D_{N\alpha}
\left(
\left(
\tf{N}{2\pi}
\right)^{(m+r)/2}t
\right).
\]
We should note that the total mass of the measure $\nu_{N,r}$ is $\vol(M_{P})(N/2\pi)^{(m+r)/2}$, and
hence it tends to infinity as $N$ goes to infinity.
But its $k$-th moment satisfies, for each positive integer $k$,
\begin{equation}
\int x^{k}\,d\nu_{N,k}(x)=
\frac{c(P,\alpha)^{k-1}}{k^{(m+r)/2}}
(1+O_{k}(N^{-1}))
\to \frac{c(P,\alpha)^{k-1}}{k^{(m+r)/2}}\quad (N\to \infty).
\label{mom1}
\end{equation}

By Propositions \ref{limitS} and \ref{limitV},
the support of the measure $\nu_{N,r}$ is contained
in a bounded interval in $[0,+\infty)$ which is independent of $N$.

Furthermore, by \eqref{mom1}, the measure $xd\nu_{N,r}(x)$
is a finite measure on the real line whose support lies in
a bounded interval in $[0,+\infty)$ which is independent of $N$.
This implies that the sequence of the finite measures $xd\nu_{N,r}(x)$ on the real
line has a weak limit, say $\mu$, and it must satisfy

\begin{equation}
\int x^{k}\,d\mu(x)=
\frac{c(P,\alpha)^{k}}{(k+1)^{(m+r)/2}},\quad k=0,1,2,\ldots.
\label{mom2}
\end{equation}

Thus, the weak limits of the measures $xd\nu_{N,r}(x)$ are probability measures,
and they have supports contained in the interval $[0,c(P,\alpha)]$.

\begin{lem}
Let $c$ be a positive constant, and let $h$ be a positive integer.
Let $\mu$ be a probability measure on the real line such that
\[
\int x^{k}\,d\mu(x)=\frac{c^{k}}{(k+1)^{h/2}},\quad k \geq 0.
\]
Then $\mu$ must coincide with the measure $\rho_{c,h}(x)\,dx$ where
\begin{equation}
\rho_{c,h}(x)=\frac{1}{c\Gamma (h/2)}\chi_{(0,c)}(x)(\log (c/x))^{h/2-1}.
\label{density1}
\end{equation}
\label{momentlemma}
\end{lem}

\begin{proof}
We shall use the following formula:
\[
\frac{1}{w^{s}}=\frac{1}{\Gamma (s)}\int_{0}^{\infty}e^{-wt}t^{s-1}\,dt
\]
to compute the Fourier transform
\[
\hat{\mu}(\xi)=\int e^{ix\xi}\,d\mu(x)
\]
of the measure $\mu$ satisfying the condition in the lemma.
Substituting $s=h/2$, $w=k$ in the above formula, we get
\[
\begin{split}
\hat{\mu}(\xi)=\sum_{k \geq 0}\frac{(ic\xi)^{k}}{k!(k+1)^{h/2}}
= & \frac{1}{\Gamma (h/2)}
\int_{0}^{\infty} e^{ic\xi e^{-t}}e^{-t}t^{h/2-1}\,dt \\
= & \frac{1}{c\Gamma (h/2)}
\int_{0}^{c}e^{i\xi x}(\log (c/x))^{h/2-1}\,dx,
\end{split}
\]
and hence $\hat{\mu}(\xi)=\hat{\rho_{c,h}}(\xi)$. This completes the proof.
\end{proof}

As a corollary to Lemma \ref{momentlemma}, we have the following.

\begin{cor}
The sequence of the finite measures $\{xd\nu_{N,r}(x) \}$
converges weakly to the probability measure $\rho_{c(P,\alpha),r}\,dx$, where
the density $\rho_{c(P,\alpha),r}$ is given by \eqref{density1} with $c=c(P,\alpha)$
and $h=m+r$.
\label{Cm}
\end{cor}

\vspace{3pt}

\noindent{{\it Proof of Theorem \ref{mldist}.}} \hspace{3pt}
The rescaled distribution function $F_{N\alpha}^{r}(t)$ for $t>0$ of
the measure $\nu_{N}=|f_{N\alpha}|^{2}_{*}{\rm dv}^{r}_{N}$ is given by
\[
F_{N\alpha}^{r}(t)=\int \chi_{(t,+\infty)}(x)\,d\nu_{N,r}(x),
\]
where $\chi_{(t,+\infty)}$ is the characteristic function of the interval
$(t,+\infty)$. Now set $d\mu_{N,r}(x)=xd\nu_{N,r}(x)$ and write
\[
F_{N\alpha}^{r}(t) =
\int x^{-1}\chi_{(t,+\infty)}(x)\,d\mu_{N,r}(x).
\]
By Corollary \ref{Cm}, we know that the sequence of finite measures $\mu_{N}$ converges
weakly to the probability measure $\rho_{c(P,\alpha),r}(x)\,dx$.
For fixed $t>0$,
by approximating the function $x^{-1}\chi_{(t,+\infty)}(x)$
by a sequence of continuous functions, and by using the Lebesgue
convergence theorem, one easily obtains that, for every $0<t\leq c(P,\alpha)$,
\[
\begin{split}
\lim_{N \to +\infty}F_{N\alpha}^{r}(t) &=\int x^{-1}\chi_{(t,+\infty)}(x)\rho_{c(P,\alpha),r}(x)\,dx \\
&=\frac{1}{c(P,\alpha)\Gamma ((m+r)/2)}\int_{t}^{c(P,\alpha)}x^{-1}(\log (c(P,\alpha)/x))^{(m+r)/2-1}\,dx \\
&=\frac{1}{c(P,\alpha)\Gamma ((m+r)/2)}\int_{0}^{\log (c(P,\alpha)/t)}s^{(m+r)/2-1}\,ds \\
&=\frac{1}{c(P,\alpha)\Gamma ((m+r)/2 +1)}(\log (c(P,\alpha)/t))^{(m+r)/2}.
\end{split}
\]
\hfill\qedsymbol

\subsection{Non-rescaled distribution functions}

In the previous section, we derived the limit of the rescaled distribution functions
$F_{N\alpha}^{r}(t)$. In the definition of the rescaled function $F_{N\alpha}^{r}(t)$,
the parameter $t$ is rescaled by the factor $(N/2\pi)^{(m+r)/2}$
so that it tends to $+\infty$ as $N \to \infty$. Then the corresponding volume
\begin{equation}
\vol_{M_{P}} (z\,;\,|\wh\phi^{P}_{N\alpha}(z)|^2 > (N/2\pi)^{(m+r)/2}t)
\label{sample1}
\end{equation}
is of order $N^{-(m+r)/2}$. This was the reason why we need to
multiply the volume \eqref{sample1} by the extra factor
$(N/2\pi)^{(m+r)/2}$ in the definition of the rescaled
distribution functions. As a result, the limit distribution has a
universal form (Theorem \ref{mldist}). However, one may ask, of
course, what the limit of the non-rescaled distributions
$D_{N\alpha}(t)$ is. But, for each fixed $t >0$, $D_{N\alpha}(t)
\to 0$ as $N \to \infty$ by Theorem \ref{pointww} and Proposition
\ref{LOCALSS}. Therefore, we need to replace $D_{N\alpha}(t)$ by
$D_{N\alpha}(W_{N}(t))$ for an appropriate sequence $\{W_{N}(t)\}$
of positive functions in $t>0$ which compensate for the flattening
rate. Our sequence will satisfy the conditions $W_{N}(t) \to 0$ as
$N \to \infty$ and
\begin{equation}W_{N}(t)
\to 0,\quad
W_{N}(t)^{-1/N} \to W(t) \quad (N \to \infty), \label{weight}
\end{equation}
for a function $W(t) >1$. We then  have the following limit
distribution law:

\begin{theo}\label{thnonrescaled}
Let $\{W_N\}$ satisfy \eqref{weight}.  Then
\begin{equation}
\lim_{N \to \infty}D_{N}(W_{N}(t))= \vol_{M_{P}} ((\xi,\rho) \in
\C^{r} \times (\R^{*})^{m-r}\,;\,\Psi_{\alpha}(\xi,\rho) < \log W(t) ).
\label{nonrescaled}
\end{equation}
In particular, by taking the function $W_{N}$ as $W_{N}(t)=e^{-Nt}$,
we have
\begin{equation}
\lim_{N \to \infty} D_{N}(e^{-Nt})=
\frac{1}{\pi^{r}}\int_{\{(\xi, \rho) \in
\C^{r} \times (\R^{*})^{m-r}\,;\,\Psi_{\alpha}^{P}(\xi,\rho) < t \}} L(\xi,\rho)\,dm(\xi)d\rho.
\label{expfunc3}
\end{equation}
If $\alpha$ is in the interior $P^{o}$ of the polytope $P$, we have
\begin{equation}
\lim_{N\to \infty}
D_{N}(e^{-Nt})=
\int_{\{\rho \in \R^{m}\,;\,b_{\alpha}^{P}(\rho) <t\}}
\det A(\rho)\,d\rho.
\label{INTnon}
\end{equation}
\label{tnonrescaled}

\end{theo}

\begin{proof}
The proof of \eqref{INTnon} is the same as that for \eqref{nonrescaled} and \eqref{expfunc3},
and \eqref{expfunc3} follows from \eqref{nonrescaled} and Lemma \ref{volume22}.
Thus, we give a proof only for \eqref{nonrescaled}.
First of all, we note that the set $U_{v_{0}}$ introduced in Section \ref{massasym} is dense in $M_{P}$,
and hence we may consider the volume of the set
\[
S_{N}(t):=\{(\xi,\zeta) \in \C^{r} \times (\C^{*})^{m-r}\,;\,
|\wh\phi^{P}_{N\alpha}(\xi,\zeta)|^{2} >W_N(t)\},
\]
so that $D_{N\alpha}(W_N(t))=\vol_{M_{P}} (S_{N}(t))$.
Let $(\xi,\zeta) \in S_{N}$.
We write $\zeta =e^{\rho/2 +i\theta}$.
Then, by Proposition \ref{LOCALSS}, we have
\begin{equation}
\label{estimate00}
\Psi_{\alpha}(\xi,\rho) < \log
\left(
CN^{-(m+r)/2}W_{N}(t)(1+a_{N}(\xi,\zeta))
\right)^{-1/N},
\end{equation}
where $C$ is a constant, and $a_{N}$ is a function of order $O(N^{-1})$.
Thus, we obtain
\[
\begin{split}
D_{N\alpha}&(W_{N}(t))\\
=&\vol_{M_{P}}
\left(
(\xi,\zeta) \in \C^{r} \times (\C^{*})^{m-r}\,;\,
\Psi_{\alpha}(\xi,\rho)
<\log
\left(
CN^{-(m+r)/2}W_{N}(t)(1+a_{N}(\xi,\zeta))
\right)^{-\frac{1}{N}}
\right).
\end{split}
\]
This combined with Lemma \ref{volume22} shows the assertion.
\end{proof}

\vspace{5pt}

Our final aim is to prove theorem \ref{unrescal}, which gives the
asymptotic limit of the distribution function $D_{N}(t)$ itself
without any rescaling. Since, by Theorem \ref{pointww} and
Propositions \ref{LOCALSS} and \ref{LOCALV}, the monomial
$|\wh\phi_{N\alpha}^{P}|^{2}$ decays exponentially away from the
corresponding invariant torus $\mu_{P}^{-1}(\alpha)$, it is
obvious that $D_{N}(t)$ tends to zero as $N \to \infty$. So, a
problem is to find the decay rate of $D_{N}(t)$ for any (but
fixed) $t >0$.

For every $N$ and $0 \leq r \leq m$, we set
\[
s_{N}=\log
\left(
\tf{N}{2\pi}
\right)^{(m+r)/2}
=\tf{m+r}{2}\log
\left(
\tf{N}{2\pi}
\right), \quad
r_{N}=
\left(
\tf{s_{N}}{N}
\right)^{1/2}.
\]

\begin{lem}
\label{Lest1}
Let $t >0$. Let $\alpha \in P$ be a lattice point with $\dim \mu_{P}^{-1}(\alpha)=m-r$
and lie in a face $F$ with $\dim F=m-r$.
Let $(\xi_{N},\zeta_{N}) \in \C^{r} \times (\C^{*})^{m-r}$ satisfy
$|\wh{\varphi}_{N\alpha}^{P}(\xi_{N},\zeta_{N})|^{2}>t$.
We write $\zeta_{N}=e^{(\rho_{\alpha}^{F}+u_{N})/2+i\theta_{N}}$.
Then, we have
\[
|\xi_{N}|^{2}+|u_{N}|^{2}
=O(r_{N}^{2})=
O(N^{-1}\log N)
\]
locally uniformly in $t>0$.
\end{lem}

\begin{proof}
For every $c>0$, we set
\begin{equation}
\label{BALL}
B_{\alpha}(c)=\{
(\xi,e^{(\rho_{\alpha}^{F}+u)/2}+i\theta) \in
\C^{r} \times (\C^{*})^{m-r}\,;\,
|\xi|^{2} +|u|^{2} <c^{2}
\}.
\end{equation}
Then, by the argument in the proof of Proposition \ref{LOCALSS},
we can find positive constants $c_{0}$, $C_{0}$ such that
\[
c_{0}(|\xi|^{2}+|u|^{2}) \leq \Psi_{\alpha}(\xi,\rho_{\alpha}^{F}+u)
\leq C_{0}(|\xi|^{2}+|u|^{2}) \mbox{ on }B_{\alpha}(c).
\]
Therefore, by Proposition \ref{LOCALSS}, we have
\[
t \leq Ce^{s_{N}-N\Psi_{\alpha}(\xi_{N},\rho_{\alpha}^{F}+u_{N})}
\leq Ce^{s_{N}-c_{0}N(|\xi_{N}|^{2}+|u_{N}|^{2})}
\]
with some constant $C$.
Therefore, we obtain
\[
|\xi_{N}|^{2}+|u_{N}|^{2} \leq \frac{C}{N}
(\log(C/t)+s_{N}),
\]
which implies the assertion.
\end{proof}

Thus, to obtain the asymptotic estimate of the distribution function $D_{N}(t)$,
we need to find that of the monomial $|\wh\phi_{N\alpha}^{P}|^{2}$ on the ball $B_{\alpha}(cr_{N})$
of radius $O(r_{N})$ around the point $(0,\rho_{\alpha}^{F})$, where $B_{\alpha}(cr_{N})$ is defined in
\eqref{BALL}.

\begin{lem}
Let $c>0$. We denote points $B_{\alpha}(cr_{N})$ as
$(\xi,\zeta)=(r_{N}w,e^{(\rho_{\alpha}^{F}+r_{N}u)/2+i\theta})$
with $|w|^{2}+|u|^{2} \leq c^{2}$.
Then we have
\begin{equation}
\begin{split}
|\wh{\varphi}_{N\alpha}^{P}(r_{N}w,&\rho_{\alpha}^{F}+r_{N}u)|^{2}\\
=&c(P,\alpha)e^{s_{N}-s_{N}\ispa{H\Psi_{\alpha}(0,\rho_{\alpha}^{F})(w,u),(w,u)}/2}
(1+O(N^{-1/2}(\log N)^{3/2})),
\end{split}
\end{equation}
where the Hessian $H\Psi_{\alpha}(0,\rho_{\alpha}^{F})$ of the function $\Psi_{\alpha}$
on $\C^{r} \times \R^{m-r}$
at the point $(0,\rho_{\alpha}^{F})$ is given by \eqref{hess222}
\label{Logest}
\end{lem}

\begin{proof}
By a Taylor expansion, we have
\[
\Psi_{\alpha}(r_{N}w,\rho_{\alpha}^{F}+r_{N}u)=
\frac{r_{N}^{2}}{2}
\ispa{H\Psi_{\alpha}(0,\rho_{\alpha}^{F})(w,u),(w,u)}
+R(r_{N}(w,u))
\]
for $|w|^{2}+|u|^{2} \leq c^{2}$,
where $R(r_{N}(w,u))=O(r_{N}^{3})$.
In particular, we have
\[
\begin{split}
e^{-N\Psi_{\alpha}(r_{N}w,\rho_{\alpha}^{F}+r_{N}u)}
&=e^{-s_{N}\ispa{H\Psi_{\alpha}(0,\rho_{\alpha}^{F})(w,u),(w,u)}/2}
(1+O(Nr_{N}^{3})) \\
&=e^{-s_{N}\ispa{H\Psi_{\alpha}(0,\rho_{\alpha}^{F})(w,u),(w,u)}/2}
(1+O(N^{-1/2}(\log N)^{3/2})).
\end{split}
\]
>From this and the estimate in Proposition \ref{LOCALSS}, the
assertion follows for $1 \leq r \leq m-1$. For $r=0,m$, precisely
the same argument replacing $\Psi_{\alpha}$ by $b_{\alpha}^{P}$ or
$\log K$ with Theorem \ref{pointww} and Proposition \ref{LOCALV}
will show the assertion.
\end{proof}

\noindent{{\it Proof of Theorem \ref{unrescal}.}}\hspace{5pt}
For every $t>0$, we set
\[
S_{N}(t)=\{
\eta=(\xi,\zeta) \in \C^{r} \times (\C^{*})^{m-r}\,;\,
|\wh{\varphi}_{N\alpha}^{P}(\xi,\zeta)|^{2} >t
\}.
\]
Then, by Lemma \ref{Lest1},
there is a constant $c>0$ (depending on a fixed $t >0$) such that
$S_{N}(t) \subset B_{\alpha}(cr_{N})$.
Let $(\xi,\zeta)=(r_{N}w,e^{(\rho_{\alpha}^{F}+r_{N}u})/2 +i\theta) \in B_{\alpha}(cr_{N})$.
Then, by Lemma \ref{Logest},
$(\xi,\zeta) \in S_{N}(t)$ if and only if
\[
t < c(P,\alpha)e^{s_{N}-s_{N}\ispa{H\Psi_{\alpha}(0,\rho_{\alpha}^{F})(w,u),(w,u)}/2}(1+a_{N}(w,u)),
\]
where $a_{N}(w,u)$ is a function of order $N^{-1/2}(\log N)^{3/2}$ uniformly in
$(w,u)$ with $|w|^{2} +|u|^{2} \leq c^{2}$,
and $c(P,\alpha)$ is define in \eqref{CONSTC}.
This is equivalent to the following estimate:
\[
\ispa{H\Psi_{\alpha}(w,u),(w,u)}/2 < 1 +\frac{1}{s_{N}}\log
\left(
\frac{c(P,\alpha)}{t}(1+a_{N}(w,u))
\right).
\]
Therefore, by Lemma \ref{volume22}, we have
\[
D_{N\alpha}(t)=
\frac{1}{\pi^{r}}
\int_{\{
(\xi,\rho)=(r_{N}w,\rho_{\alpha}^{P}+r_{N}u)\,;\,
\ispa{H(w,u),(w,u)}/2
<1 +\frac{1}{s_{N}}
\log(c(P,\alpha)(1+a(w,u))/t)
\}}
L(\xi,\rho)\,dm(\xi)d\rho,
\]
where we set, for simplicity, $H=H\Psi_{\alpha}(0,\rho_{\alpha}^{F})$.
Changing the variables $(\xi,\rho)$ to $(w,u)$, we have
\[
D_{N\alpha}(t)
=\frac{r_{N}^{m+r}}{\pi^{r}}
\int_{\{
(w,u) \,;\,
\ispa{H(w,u),(w,u)}/2 <
1+ \frac{1}{s_{N}}
\log(c(P,\alpha)(1+a_{N}(w,u))/t)\}}
\!\!\!\!\!\!L(r_{N}w,\rho_{\alpha}^{F}+r_{N}u)\,dm(w)du.
\]
A Taylor expansion gives
\[
L(r_{N}w,\rho_{\alpha}^{F}+r_{N}u)=
L(0,\rho_{\alpha}^{F})+R_{N}(w,u)
\]
with the error term $R_{N}(w,u)=O(r_{N})$ uniformly
for $|w|^{2}+|u|^{2} \leq c^{2}$.
Thus, we obtain
\[
r_{N}^{-(m+r)}D_{N\alpha}(t)=
\frac{L(0,\rho_{\alpha}^{F})}{\pi^{r}}
\int_{\{
(w,u) \,;\,
\ispa{H(w,u),(w,u)}/2 <
1+ \frac{1}{s_{N}}
\log(c(P,\alpha)(1+a_{N}(w,u))/t)\}}
\!\!\!\!\!dm(w)du +O(r_{N}).
\]

Note that $s_{N} \to \infty$ while $r_{N} \to 0$ as $N \to \infty$.
The function $a_{N}(w,u)$ tends to zero uniformly in $|w|^{2}+|u|^{2} \leq c^{2}$.
Therefore, we conclude, for a fixed $t>0$,
\[
\begin{split}
\lim_{N \to \infty}r_{N}^{-(m+r)}D_{N\alpha}(t)&=
\frac{L(0,\rho_{\alpha}^{F})}{\pi^{r}}
\int_{\{(w,u) \in \C^{r} \times \R^{m-r}\,;\,
\ispa{H(w,u),(w,u)}/2 <1\}}
dm(w)du. \\
&=\frac{2^{(m+r)/2}}{\pi^{r}}
\frac{L(0,\rho_{\alpha}^{F})}{\sqrt{\det H}}
\vol_{M_{P}} (x \in \R^{m+r}\,;\,|x| <1).
\end{split}
\]
By \eqref{DET} and the well-known formula for the volume of the unit disk in $\R^{m+r}$,
we conclude the assertion.
\hfill\qedsymbol

\end{document}